\documentclass[11pt]{article}
\usepackage{amsmath, upgreek}
\usepackage{amssymb}
\usepackage{amsfonts}
\usepackage{amsmath,tabu}
\usepackage{cite}
\usepackage{color}
\usepackage[final]{pdfpages}
\usepackage{xcolor} 
\usepackage{amscd}
\usepackage{xspace}
\usepackage{verbatim}
\usepackage{graphicx}
\newcommand{\mysection}[1]{
\section{#1}\setcounter{equation}{0}}
\title{\bf Singularities and asymptotics of solutions of the\\
Chandrasekhar-Hamilton-Jacobi equation
}
\author{{\bf Marie-Fran\c{c}oise Bidaut-V\'eron\footnote{Email: veronmf@univ-tours.fr}, Laurent V\'eron\footnote{Corresponding author. Email: veronl@univ-tours.fr}}}
 
\date{}
\begin{document}

 \maketitle


\newcommand{\txt}[1]{\;\text{ #1 }\;}
\newcommand{\tbf}{\textbf}
\newcommand{\tit}{\textit}
\newcommand{\tsc}{\textsc}
\newcommand{\trm}{\textrm}
\newcommand{\mbf}{\mathbf}
\newcommand{\mrm}{\mathrm}
\newcommand{\bsym}{\boldsymbol}
\newcommand{\scs}{\scriptstyle}
\newcommand{\sss}{\scriptscriptstyle}
\newcommand{\txts}{\textstyle}
\newcommand{\dsps}{\displaystyle}
\newcommand{\fnz}{\footnotesize}
\newcommand{\scz}{\scriptsize}
\newcommand{\be}{\begin{equation}}
\newcommand{\bel}[1]{\begin{equation}\label{#1}}
\newcommand{\ee}{\end{equation}}
\newcommand{\eqnl}[2]{\begin{equation}\label{#1}{#2}\end{equation}}
\newcommand{\barr}{\begin{eqnarray}}
\newcommand{\earr}{\end{eqnarray}}
\newcommand{\bars}{\begin{eqnarray*}}
\newcommand{\ears}{\end{eqnarray*}}
\newcommand{\nnu}{\nonumber \\}
\newtheorem{subn}{\name}
\renewcommand{\thesubn}{}
\newcommand{\bsn}[1]{\def\name{#1}\begin{subn}}
\newcommand{\esn}{\end{subn}}
\newtheorem{sub}{\name}[section]
\newcommand{\dn}[1]{\def\name{#1}}   
\newcommand{\bs}{\begin{sub}}
\newcommand{\es}{\end{sub}}
\newcommand{\bsl}[1]{\begin{sub}\label{#1}}
\newcommand{\bth}[1]{\def\name{Theorem}
\begin{sub}\label{t:#1}}
\newcommand{\blemma}[1]{\def\name{Lemma}
\begin{sub}\label{l:#1}}
\newcommand{\bcor}[1]{\def\name{Corollary}
\begin{sub}\label{c:#1}}
\newcommand{\bdef}[1]{\def\name{Definition}
\begin{sub}\label{d:#1}}
\newcommand{\bprop}[1]{\def\name{Proposition}
\begin{sub}\label{p:#1}}

\newcommand{\R}{\eqref}
\newcommand{\rth}[1]{Theorem~\ref{t:#1}}
\newcommand{\rlemma}[1]{Lemma~\ref{l:#1}}
\newcommand{\rcor}[1]{Corollary~\ref{c:#1}}
\newcommand{\rdef}[1]{Definition~\ref{d:#1}}
\newcommand{\rprop}[1]{Proposition~\ref{p:#1}}
\newcommand{\BA}{\begin{array}}
\newcommand{\EA}{\end{array}}
\newcommand{\BAN}{\renewcommand{\arraystretch}{1.2}
\setlength{\arraycolsep}{2pt}\begin{array}}
\newcommand{\BAV}[2]{\renewcommand{\arraystretch}{#1}
\setlength{\arraycolsep}{#2}\begin{array}}
\newcommand{\BSA}{\begin{subarray}}
\newcommand{\ESA}{\end{subarray}}
\newcommand{\BAL}{\begin{aligned}}
\newcommand{\EAL}{\end{aligned}}
\newcommand{\BALG}{\begin{alignat}}
\newcommand{\EALG}{\end{alignat}}
\newcommand{\BALGN}{\begin{alignat*}}
\newcommand{\EALGN}{\end{alignat*}}
\newcommand{\note}[1]{\textit{#1.}\hspace{2mm}}
\newcommand{\Proof}{\note{Proof}}
\newcommand{\qeda}{\hspace{10mm}\hfill $\square$}
\newcommand{\qed}{\\
${}$ \hfill $\square$}
\newcommand{\Remark}{\note{Remark}}
\newcommand{\modin}{$\,$\\[-4mm] \indent}
\newcommand{\forevery}{\quad \forall}
\newcommand{\set}[1]{\{#1\}}
\newcommand{\setdef}[2]{\{\,#1:\,#2\,\}}
\newcommand{\setm}[2]{\{\,#1\mid #2\,\}}
\newcommand{\mt}{\mapsto}
\newcommand{\lra}{\longrightarrow}
\newcommand{\lla}{\longleftarrow}
\newcommand{\llra}{\longleftrightarrow}
\newcommand{\Lra}{\Longrightarrow}
\newcommand{\Lla}{\Longleftarrow}
\newcommand{\Llra}{\Longleftrightarrow}
\newcommand{\warrow}{\rightharpoonup}
\newcommand{
\paran}[1]{\left (#1 \right )}
\newcommand{\sqbr}[1]{\left [#1 \right ]}
\newcommand{\curlybr}[1]{\left \{#1 \right \}}
\newcommand{\abs}[1]{\left |#1\right |}
\newcommand{\norm}[1]{\left \|#1\right \|}
\newcommand{
\paranb}[1]{\big (#1 \big )}
\newcommand{\lsqbrb}[1]{\big [#1 \big ]}
\newcommand{\lcurlybrb}[1]{\big \{#1 \big \}}
\newcommand{\absb}[1]{\big |#1\big |}
\newcommand{\normb}[1]{\big \|#1\big \|}
\newcommand{
\paranB}[1]{\Big (#1 \Big )}
\newcommand{\absB}[1]{\Big |#1\Big |}
\newcommand{\normB}[1]{\Big \|#1\Big \|}
\newcommand{\produal}[1]{\langle #1 \rangle}

\newcommand{\thkl}{\rule[-.5mm]{.3mm}{3mm}}
\newcommand{\thknorm}[1]{\thkl #1 \thkl\,}
\newcommand{\trinorm}[1]{|\!|\!| #1 |\!|\!|\,}
\newcommand{\bang}[1]{\langle #1 \rangle}
\def\angb<#1>{\langle #1 \rangle}
\newcommand{\vstrut}[1]{\rule{0mm}{#1}}
\newcommand{\rec}[1]{\frac{1}{#1}}
\newcommand{\opname}[1]{\mbox{\rm #1}\,}
\newcommand{\supp}{\opname{supp}}
\newcommand{\dist}{\opname{dist}}
\newcommand{\myfrac}[2]{{\displaystyle \frac{#1}{#2} }}
\newcommand{\myint}[2]{{\displaystyle \int_{#1}^{#2}}}
\newcommand{\mysum}[2]{{\displaystyle \sum_{#1}^{#2}}}
\newcommand {\dint}{{\displaystyle \myint\!\!\myint}}
\newcommand{\q}{\quad}
\newcommand{\qq}{\qquad}
\newcommand{\hsp}[1]{\hspace{#1mm}}
\newcommand{\vsp}[1]{\vspace{#1mm}}
\newcommand{\ity}{\infty}
\newcommand{\prt}{\partial}
\newcommand{\sms}{\setminus}
\newcommand{\ems}{\emptyset}
\newcommand{\ti}{\times}
\newcommand{\pr}{^\prime}
\newcommand{\ppr}{^{\prime\prime}}
\newcommand{\tl}{\tilde}
\newcommand{\sbs}{\subset}
\newcommand{\sbeq}{\subseteq}
\newcommand{\nind}{\noindent}
\newcommand{\ind}{\indent}
\newcommand{\ovl}{\overline}
\newcommand{\unl}{\underline}
\newcommand{\nin}{\not\in}
\newcommand{\pfrac}[2]{\genfrac{(}{)}{}{}{#1}{#2}}

\def\ga{\alpha}     \def\gb{\beta}       \def\gg{\gamma}
\def\gc{\chi}       \def\gd{\delta}      \def\ge{\epsilon}
\def\gth{\theta}                         \def\vge{\varepsilon}
\def\gf{\phi}       \def\vgf{\varphi}    \def\gh{\eta}
\def\gi{\iota}      \def\gk{\kappa}      \def\gl{\lambda}
\def\gm{\mu}        \def\gn{\nu}         \def\gp{\pi}
\def\vgp{\varpi}    \def\gr{\rho}        \def\vgr{\varrho}
\def\gs{\sigma}     \def\vgs{\varsigma}  \def\gt{\tau}
\def\gu{\upsilon}   \def\gv{\vartheta}   \def\gw{\omega}
\def\gx{\xi}        \def\gy{\psi}        \def\gz{\zeta}
\def\Gg{\Gamma}     \def\Gd{\Delta}      \def\Gf{\Phi}
\def\Gth{\Theta}
\def\Gl{\Lambda}    \def\Gs{\Sigma}      \def\Gp{\Pi}
\def\Gw{\Omega}     \def\Gx{\Xi}         \def\Gy{\Psi}

\def\CS{{\mathcal S}}   \def\CM{{\mathcal M}}   \def\CN{{\mathcal N}}
\def\CR{{\mathcal R}}   \def\CO{{\mathcal O}}   \def\CP{{\mathcal P}}
\def\CA{{\mathcal A}}   \def\CB{{\mathcal B}}   \def\CC{{\mathcal C}}
\def\CD{{\mathcal D}}   \def\CE{{\mathcal E}}   \def\CF{{\mathcal F}}
\def\CG{{\mathcal G}}   \def\CH{{\mathcal H}}   \def\CI{{\mathcal I}}
\def\CJ{{\mathcal J}}   \def\CK{{\mathcal K}}   \def\CL{{\mathcal L}}
\def\CT{{\mathcal T}}   \def\CU{{\mathcal U}}   \def\CV{{\mathcal V}}
\def\CZ{{\mathcal Z}}   \def\CX{{\mathcal X}}   \def\CY{{\mathcal Y}}
\def\CW{{\mathcal W}} \def\CQ{{\mathcal Q}}
\def\BBA {\mathbb A}   \def\BBb {\mathbb B}    \def\BBC {\mathbb C}
\def\BBD {\mathbb D}   \def\BBE {\mathbb E}    \def\BBF {\mathbb F}
\def\BBG {\mathbb G}   \def\BBH {\mathbb H}    \def\BBI {\mathbb I}
\def\BBJ {\mathbb J}   \def\BBK {\mathbb K}    \def\BBL {\mathbb L}
\def\BBM {\mathbb M}   \def\BBN {\mathbb N}    \def\BBO {\mathbb O}
\def\BBP {\mathbb P}   \def\BBR {\mathbb R}    \def\BBS {\mathbb S}
\def\BBT {\mathbb T}   \def\BBU {\mathbb U}    \def\BBV {\mathbb V}
\def\BBW {\mathbb W}   \def\BBX {\mathbb X}    \def\BBY {\mathbb Y}
\def\BBZ {\mathbb Z}

\def\GTA {\mathfrak A}   \def\GTB {\mathfrak B}    \def\GTC {\mathfrak C}
\def\GTD {\mathfrak D}   \def\GTE {\mathfrak E}    \def\GTF {\mathfrak F}
\def\GTG {\mathfrak G}   \def\GTH {\mathfrak H}    \def\GTI {\mathfrak I}
\def\GTJ {\mathfrak J}   \def\GTK {\mathfrak K}    \def\GTL {\mathfrak L}
\def\GTM {\mathfrak M}   \def\GTN {\mathfrak N}    \def\GTO {\mathfrak O}
\def\GTP {\mathfrak P}   \def\GTR {\mathfrak R}    \def\GTS {\mathfrak S}
\def\GTT {\mathfrak T}   \def\GTU {\mathfrak U}    \def\GTV {\mathfrak V}
\def\GTW {\mathfrak W}   \def\GTX {\mathfrak X}    \def\GTY {\mathfrak Y}
\def\GTZ {\mathfrak Z}   \def\GTQ {\mathfrak Q}

\font\Sym= msam10 
\def\SYM#1{\hbox{\Sym #1}}
\newcommand{\bdw}{\prt\Gw\xspace}
\maketitle\medskip

\noindent{\small {\bf Abstract} We study the local properties of positive solutions of the equation $-\Gd u+ m\abs{\nabla u}^q-e^{u}=0$ in a punctured domain $\Gw\setminus\{0\}$ of $\BBR^N$, $N\geq 2$, where $m$ is a positive parameter  and $q>1$. We study particularly the local behaviour of solutions with an isolated singularity or the asymptotic behaviour for solutions defined in an exterior domain. These behaviours change drastically according $q$ is smaller or larger than $2$.
}\medskip

\tableofcontents
\vspace{1mm}
\hspace{.05in}
\medskip
\mysection{Introduction}
The aim of this article is to study the behaviour near an isolated singularities of functions which satisfy the viscous {\it Chandrasekhar-Hamilton-Jacobi}
\bel{Na-1}
-\Gd u+ m|\nabla u|^q-e^{u}=0
\ee
in a punctured domain of $\BBR^N$ or in an exterior domain. Throughout this article we assume that $m>0$ and $q>1$ However since our results depend strongly on the value of $q$ with respect to $\frac N{N-1}$ and $2$, in most statement we will recall the range of the parameter $q$. Since the problem is local and invariant by translation, we assume that $(\ref{Na-1})$ holds
in $B_\gr\setminus\{0\}$. This equation exhibits a large variety of phenomena, since besides the diffusion operator,  the reaction terms are nonlinear and of a different nature. According to the range of values of the parameters and the solutions, the behaviour can be modelled by three underlying equations:\\
the {\it Ricatti type equation or viscous Hamilton-Jacobi equation} 
\bel{Na-2}
-\Gd u+ m|\nabla u|^q=0,
\ee
the {\it eikonal equation}
\bel{Na-3}
m|\nabla u|^q-e^{u}=0,
\ee
and the {\it Emden equation or Chandrasekhar equation} 
\bel{Na-4}
-\Gd u-e^{u}=0.
\ee
But at each occurrence the difficulty is to prove that one of the three terms in the equation is negligible compared to the two others. Equation $(\ref{Na-1})$ is a model for more general equations like
\bel{Na-5}
-\Gd u+ f(|\nabla u|)-g(u)=0
\ee
where $f$ and $g$ are positive continuous superlinear functions. The specific effects we aim to put into light need explicit model functions such as the one in  $(\ref{Na-1})$. The guiding thread of the equations $(\ref{Na-2})$, $(\ref{Na-3})$, $(\ref{Na-4})$ is the existence of specific  solutions invariant under some scaling transformation specific to each equation.  \\
For the Hamilton-Jacobi equation the transformation is $T^R_\ell$ defined for $\ell>0$ by 
\bel{Na-6}
T^R_\ell[u](x)=\ell^\gb u(\ell x) \quad\text{where }\gb=\frac{2-q}{q-1}.
\ee
Self-similar solutions of $(\ref{Na-2})$ are expressed in spherical coordinates $(r,\gs)\in \BBR_+\ti S^{N-1}$ by 
\bel{Na-7}
u(r,\gs)=r^{-\gb}\gw(\gs)
\ee
where $\gw$ satisfies 
\bel{Na-8}
-\Gd_g\gw-\gb(\gb+2-N)\gw-m(\gb^2\gw^2+|\nabla_g\gw|^2)^\frac q2=0\quad\text{in } S^{N-1}
\ee
where $\Gd_g$ and $\nabla_g$ are respectively the Laplace-Beltrami operator and the covariant gradient on the unit sphere $S^{N-1}$ endowed with the canonical metric $g$ induced by the Euclidean metric thanks to the isometrical imbedding of $S^{N-1}$ into $\BBR^N$. \\
For the eikonal equation the transformation is $T^{ei}_\ell$ defined for $\ell>0$ by 
\bel{Na-9}
T^{ei}_\ell[u](x)=q\ln \ell +u(\ell x).
\ee
Self-similar solutions of $(\ref{Na-3})$ are expressed by 
\bel{Na-10}
u(r,\gs)=-q\ln r+\gw(\gs)
\ee
where $\gw$ satisfies 
\bel{Na-11}
m(q^2+|\nabla_g\gw|^2)^\frac q2-e^\gw=0\quad\text{in } S^{N-1}.
\ee
It is clear  that the only solution is the constant $\gw_e=\ln m+q\ln q$. \\
Finally, for the Chandrasekhar-Emden equation the similarity transformation is 
\bel{Na-12}
T^{E}_\ell[u](x)= u(\ell x)+2\ln\ell.
\ee
Self-similar equation of $(\ref{Na-4})$ are under the form
\bel{Na-13}
u(r,\gs)=-2\ln r+\gw(\gs) 
\ee
where $\gw$ satisfies
\bel{Na-14}
\Gd_g\gw-2(N-2)+e^\gw=0\quad\text{in } S^{N-1}.
\ee
If $N=2$ there is no solution. If $N=3$ there is an unbounded set of solutions under the form 
$\gw=\ln|det\phi|$ where $\phi$ is a conformal transformation of $S^2$ (see \cite{CY}, \cite{On}). If $4\leq N\leq 10$, besides the constant solution $\gw_E=\ln(2N-4)$, there exists an infinity of solutions 
(see \cite{Min}). \smallskip

Concerning the three above equations the question of existence of global solutions, i.e. defined in whole $\BBR^N$, is fairly understood: for the Hamilton-Jacobi equation it is proved in \cite{Lions} that any global  solution is constant. Furthermore this result still holds if the Laplacian is replaced by the $p$-Laplacian, see \cite{BVGHV-0}. Concerning the Chandrasekhar-Emden equation, it is straightforward that no global solution bounded from below exists, but infinitely many radial ones tending to  $-\infty$ at infinity; they are obtained by the shooting method. Finally concerning the eikonal equation, 
$(\ref{Na-3})$ it can be written under the more classical form
\bel{Na-3*}
|\nabla (e^{-\frac uq})|=\frac{m^{\frac 1q}}{q}.
\ee
it has been proved in \cite{CrLi} that no global bounded viscosity solution exists, or even bounded from below. Singular signed radial solutions are given in $\ref{EIK1}$.
\\

The behaviour of solutions of $(\ref{Na-1})$ differs drastically according $1<q<2$, $q=2$ and $q>2$. {\it These solutions are naturally signed as a consequence of the exponential term}. Furthermore, when $1<q<2$ the value $q=\frac N{N-1}$ is an important threshold for the Hamilton-Jacobi equation and therefore for $(\ref{Na-1})$. Note that the case $q=2$ is special in the sense that if we set $v=e^{-mu}$, then $v$ is positive and satisfies 
\bel{Na-15}
-\Gd v+mv^{1-\frac 1m}=0.
\ee
Although of a rich nature, the analysis of this equation differs from the case $q\neq 2$. When $q=2$ the behaviour of solutions depends strongly on the value of $m$, a dependence which does not hold if $q\neq 2$, since a peculiar scaling transformation could reduce the equation to the case $m=1$. We will assume $q\neq 2$ except in Section 2 in which we recall some results in this case.  For the sake of clarity we will keep the parameter $m$.
\smallskip

The difficulty in obtaining a priori estimates of singular solutions comes from the exponential term. In the case $m=0$ no such general estimates exist and  the very complete description of solutions provided in \cite{BV-V1} concerns the  solutions $u$ in $B_\gr\setminus\{0\}$ satisfying  
\bel{Na-17}
|x|^2e^u\in L^\infty(B_\gr),
\ee
an estimate which holds if $u$ is radially symmetric. \smallskip

When $q>1$, by combining  Bernstein and Keller-Osserman methods,  we first obtain a general priori estimate of the gradient depending on bounds of $e^u$. {\it Concerning the local regularity of solutions, if 
$\nabla u$ is locally bounded they are $C^2$ by standard regularity theory (and actually more)}. Therefore any solution will be assumed to be at least $C^2$ in the domain where it is considered. Concerning supersolutions or subsolutions we will assume also throughout this article that they are $C^1$ with $\Gd u$ in the sense of distributions locally integrable, or $C^2$ when necessary. 

\bth{Th2}  Let $N\geq 1$ and $q>1$. If $u$ is any  solution of $(\ref{Na-1})$ in $B_\gr(x)$, then  there exist constants $c_j=c_j(N,q,m)>0$, $j=1,2,3$ such that 
\bel{Na-2-1}
|\nabla u(x)|\leq c_1\gr^{-\frac{1}{q-1}}+c_2\max_{z\in B_{\gr}(x)}e^{\frac {u(z)}q}+c_3\max_{z\in B_{\gr}(x)}e^{\frac {u(z)}{2(q-1)}}.
\ee
\es

\nind As a consequence, when $1<q<2$ we obtain an optimal estimate of the gradient which extends the one of  \cite{BV-V1}. 


\bcor{Th1} Assume $N\geq 2$ and $1<q<2$. If $u$ is a solution of $(\ref{Na-1})$ in $B_\gr\setminus\{0\}$ satisfying $(\ref{Na-17})$, then there exists $C>0$ depending on 
$N,q,m$ and $\norm{|x|^2e^u}_{L^\infty(B_\gr)}$ such that 
\bel{Na-18}
|\nabla u(x)|\leq C|x|^{-\frac 1{q-1}}\quad\text{for all }x\in B_{\frac \gr2}\setminus\{0\}.
\ee
Therefore $u$ satisfies 
\bel{Na-19}
|u(x)|\leq \frac C\gb\left(|x|^{-\gb}-2^\gb\gr^{-\gb}\right)+\sup_{z\in\prt B_{\frac \gr2}}|u(z)|\quad\text{for all }x\in B_{\frac \gr2}\setminus\{0\},
\ee
where $\gb$ is defined in $(\ref{Na-6})$.
\es 
Note that  the assumption $(\ref{Na-17})$ is useless in our proof when $u$ is radial. \\

An important consequence of $(\ref{Na-18})$ is Harnack inequality for solutions keeping a constant sign up to an additive constant.

\bcor{Harn1} Assume $N\geq 2$, $1<q<2$ and $u$ is a solution of $(\ref{Na-1})$ in $B_\gr\setminus\{0\}$.\\
(i) If $u$  is bounded from below by some $a\in \BBR$ and satisfies $(\ref{Na-17})$, then for any $b<a$ there exists $C=C(N,q,m,a,b)>0$ such that 
\bel{Na-20}
u(x)-b\leq C\left(u(y)-b\right)\quad\text{for all }x,y\in B_{\frac\gr2}\setminus\{0\}\,\text{ such that }|x|=|y|.
\ee
(ii) If $u$ is bounded from above by some $a\in \BBR$, then for any $b>a$ $(\ref{Na-20})$ holds.
\es

When $q>2$ the estimate $(\ref{Na-18})$ holds if we assume that $|x|^\frac{q}{q-1}e^u\in L^\infty(B_\gr)$ and by integration of $(\ref{Na-18})$ we deduce that  $u$ is bounded.\medskip

The next a priori estimate, valid when $q>2$ presents the remarkable property that {\it  it is optimal and holds for any solution}. By opposition to the case $1<q<2$ (or even if $m=0$) where 
a maximal growth estimate is needed to have an a priori estimate, in that case the presence of a super-quadratic gradient term plays the role of a strong regularising absorption. This is a surprising result in view of the results concerning the Chandrasekhar equation $(\ref{Na-4})$ for which no a priori estimate could exist in view of its instability, see \cite{Far1}, \cite{Far2} and \cite{CFRS} for a related equation. The result has also to be put in parallel with the recent paper of Quittner and Souplet \cite{QS} concerning nonlinear elliptic equation without scaling invariance as it is the case with $(\ref{Na-1})$.

\bth{Th3} Assume $N\geq 2$ and $q>2$. If $u$ is a solution of $(\ref{Na-1})$ in $B_\gr\setminus\{0\}$, then there exists $C>0$ depending on $N,q,m,\gr$, but also on $u$ such that for any $x\in B_{\frac \gr2}\setminus\{0\}$ there holds
\bel{Na-20-4}\BA{lll}\dsps
(i)\qquad\qquad\qquad\qquad &e^{u(x)}\leq C|x|^{-q},\qquad\qquad\qquad\qquad\qquad\qquad\qquad\qquad\qquad\qquad\qquad\\[2mm]
(ii)\qquad\qquad\qquad\qquad &|\nabla u(x)|\leq C|x|^{-1}.
\EA\ee
\es 
The proof of this result relies on a combination of the Bernstein method, the doubling technique and the properties of the eikonal equation in $\BBR^N$. Another surprising result deals with super-solutions. By an iterative method we prove the next result valid for any $q>1$, which applies in particular to radial solutions.

\bth{Th4} Assume $N\geq 2$ and $q>1$. If $u\in C^2(\overline B_\gr\setminus\{0\})$ satisfies
\bel{Na-21}
-\Gd u +m|\nabla u|^q\geq e^u\quad\text{in }B_{\gr}\setminus\{0\},
\ee
then there exists $C=C(N,q,m,u)$ such that 
\bel{Na-22}
e^{\gm(r)}\leq Cr^{-\max\{2,q\}}\quad\text{for all }0<r\leq\frac\gr2,
\ee
where $\dsps\gm(r)=\min_{|x|=r}u(x)$.
\es

Our main result in the case $1<q<2$ is a precise description of the behaviour of solutions of $(\ref{Na-1})$ near $0$. This behaviour depends strongly of sign of these solutions: the positive ones  behave like the singular solutions of the Chandrasekhar-Emden equations, while the negative ones like the solutions of the viscous Hamilton-Jacobi equation.

\bth{Th5} Assume $1<q<2$ and $N\geq 3$. If $u$ is a solution of $(\ref{Na-1})$ such that $|x|^2e^u$ is uniformly bounded in $B_\gr$. Then\smallskip

\nind (i) Either $u$ is smooth solution. \smallskip

\nind (ii) Either there exists some $\gw$ which is a $C^2$ solution of $(\ref{Na-14})$ such that 
\bel{Na-22*}\dsps
\lim_{r\to 0}\left(u(r,\gs)-2\ln\frac 1r\right)=\gw(\gs),
\ee 
uniformly on $S^{N-1}$.\smallskip

\nind (iii) Or,  if $1<q<\frac N{N-1}$, there exists $\gg< 0$ such that 
\bel{Na-23}\dsps
\lim_{x\to 0} |x|^{N-2}u(x)=\gg.
\ee
Furthermore $e^u$ and $|\nabla u|^q$ are locally integrable in $B_{\gr}$ and there holds
 \bel{Na-24}\dsps
-\Gd u+m|\nabla u|^q-e^u=c_N\gg\gd_0\quad\text{in }\CD'(B_\gr).
\ee
\smallskip

\nind (iv) Or, if $\frac N{N-1}<q<2$ and $u$ is bounded from above, we have that
\bel{Na-25}\dsps
-c_2|x|^{-\gb}\leq u(x)\leq -c_1|x|^{-\gb},
\ee
for some positive constants $c_1,c_2$.
\smallskip

\nind (v) Or, if $q=\frac N{N-1}$ and $u$ is bounded from above, there holds
\bel{Na-26}\dsps
-c_4|x|^{2-N}(-\ln |x|)^{1-N} \leq
u(x)\leq  -c_3|x|^{2-N}(-\ln |x|)^{1-N}, 
\ee
for some positive constants $c_3,c_4$.
\es

This behaviour can be compared to the one of solutions of $(\ref{Na-1})$ when $m=0$ described in \cite[Theorem 2.1]{BV-V1}. Besides the use of the a priori estimates  the proof of $(\ref{Na-22*})$ in \rth{Th4} relies on {\it the convergence of bounded solutions of asymptotically autonomous and analytic gradient-like systems},
a result due to Huang and Takac in \cite{HuTak}  in continuation of Simon's general theory \cite{Sim}. A second fundamental tools is 
\nind{\it  the isotropy estimate} which allows to estimate the difference between a solution $u$ and its spherical average under the form
\bel{Na-34**}|u(x)-\bar u(|x|)|\leq C|x|^\gm\quad\text{for all }x\;\text{ s.t. }0<|x|\leq \frac\gr 2,
\ee
where $\gm$ is an exponent depending on the a priori estimate on $u$, see \rprop{int-3}. This isotropy estimate is obtained using {\it the integral representation of solutions of elliptic equations in an infinite cylinder.} \smallskip

When $u$ is a radial function the statements (iv) and (v) can be made more precise: (iv) has to be replaced by  \medskip

\nind {\it \nind (iv') Or, if $\frac N{N-1}<q<2$,
\bel{Na-27}\dsps
\lim_{r\to 0}r^{\gb} u(r)=\Gl_{N,m,q}:=-\frac{1}{\gb}\left(\frac{N-2-\gb}{m}\right)^{\frac 1{q-1}},
\ee}
and (v) by \smallskip

\nind {\it 
\nind (v') Or, if $q=\frac N{N-1}$,
\bel{Na-28}\dsps
\lim_{r\to 0}r^{N-2}\left(-\ln r\right)^{N-1} u(r)=-\Gl_{N,m}:=-\frac{1}{N-2}\left(\frac{m}{N-1}\right)^{1-N}.
\ee
}

When $q>2$ we prove a classification result of singular solutions of $(\ref{Na-1})$ using the estimates of \rth{Th2} and \rth{Th3}.

\bth{Th 5*} Let $N\geq 2$ and $q>2$. If $u$ is a solution of $(\ref{Na-1})$ in $B_1\setminus\{0\}$, there holds\smallskip

\nind (i) If $\dsps\lim_{x\to 0}|x|^qe^{u(x)}=0$ then $u$ can be extended as a H\"older continuous function in $B_1$ and there holds
\bel{Na-28*}\dsps
|u(x)-u(0)|\leq C|x|^{\frac{q-2}{q-1}}\quad\text{for all }x\in B_1.
\ee
Furthermore, if $u$ is radial, then the following precision holds
\bel{Na-28**}\dsps
u_r(r)=\left(\frac{N(q-1)-q}{m(q-1)r}\right)^{\frac{1}{q-1}}(1+o(1))\quad\text{as }r\to 0.
\ee
\nind (ii) If $\dsps\liminf_{x\to 0}|x|^qe^{u(x)}>0$, then 
\bel{Na-28***}
u(x)=q\ln\frac{qm^{\frac 1q}}{|x|}+O(1)\quad\text{as }x\to 0.
\ee
Furthermore, if $u$ is radial relation $(\ref{Na-28***})$ holds under the weaker condition $\dsps\limsup_{x\to 0}r^qe^{u(r)}>0$.
\es

The methods we developed for describing singularities can be used to analyse the behaviour at infinity of solutions of $(\ref{Na-1})$ in an exterior domain 
of $\BBR^N$. The main feature is a kind of exchange of the methods concerning the cases $1<q<2$ and $q>2$ at zero and at infinity. \smallskip

We first prove that if $q>1$ the estimate on supersolutions of $(\ref{Na-1})$ in $B^c_\gr$ holds under the form
\bel{Na-29}\dsps
e^{\gm(r)}\leq Cr^{-\max\{2,q\}}\quad\text{for all }r\geq 2\gr,
\ee
where $\gm(r)=\min\{u(x):|x|=r\}$.
{\it  As a consequence there exists no global solution of $(\ref{Na-1})$ which is bounded from below.} \smallskip

Concerning solutions, the estimate differs according $1<q<2$ and $q>2$. We obtain the following estimate,

\bth{th6} Let $N\geq 2$ and $1<q<2$. If $u$ is a solution of $(\ref{Na-1})$ in $B^c_1$ satisfying $\dsps \lim_{|x|\to\infty}e^{u(x)}=0$ there holds:\smallskip
\bel{Na-30}\dsps
u(x)\leq -q\ln |x|+C\quad\text{and }\,|\nabla u(x)|\leq C|x|^{-1}\quad\text{for all }\,x\in B_{2\gr}^c,
\ee
for some $C\geq 0$.
\es
The counterpart at infinity of \rth{Th3} is
\bth{th7} Let $N\geq 3$ and $q>2$. If $u$ is a solution of $(\ref{Na-1})$ in $B^c_\gr$ such that  and $|x|^2e^{u(x)}$ is uniformly bounded   there exists a $C^2$ solution $\gw$ of $(\ref{Na-14})$ such that 
\bel{Na-31}\dsps
\lim_{r\to \infty}\left(u(r,\gs)+2\ln r\right)=\gw(\gs),
\ee
uniformly on $S^{N-1}$.
\es

\nind The three main ingredients for proving \rth{th7} are:\\
{\it The two-side estimate}
\bel{Na-32}
-C|x|^{\frac{q-2}{q-1}}\leq u(x)\leq -2\ln |x|+C'\quad\text{and }\;
|\nabla u(x)|\leq C|x|^{-\frac{1}{q-1}}\quad\text{for any }x\in B_{2\gr}^c,
\ee
for some  constants $C'$ and $C$ with $C>0$,\\
\nind {\it  Harnack inequality} in the sense that for $r\geq R$ there exists $c_j=c_j(R)>0$ such that
\bel{Na-33}0<c_2\leq \frac{u(x)}{u(y)}\leq c_1\quad\text{for all }x, y\;\text{ s.t. }|x|=|y|\geq R.
\ee
\nind{\it  The isotropy result}
\bel{Na-34}|u(x)-\bar u(|x|)|\leq C|x|^{\gn}\quad\text{for all }x\;\text{ s.t. }|x|\geq R.
\ee
for some suitable $\gn$. This estimate is proved in \rprop{th*3} in a somewhat similar way as the previous one $(\ref{Na-34**})$ except that the range of the exponent $q$ is different. \medskip

In the present paper the emphasis is put on the a priori bounds or the limit behaviour of {\it any} solutions of $(\ref{Na-1})$. In the article \cite{BV-V3} we study the existence of {\it radial} singular solutions or global solutions of $(\ref{Na-1})$ introducing several types of differential systems.
Therein the methods come from the dynamical systems and in particular the stable manifold theory.

\mysection{The case $q=2$}

We give an overview of the case $q=2$ which presents a transition between the cases $1<q<2$ and $q>2$. Here the value of $m$ is fundamental. Note that many results in that case are now classical
but this is not the case of all of them. If $u$ is a solution of
\bel{q2-1}
-\Gd u+m|\nabla u|^2=e^u
\ee
in $\Gw$, we set 
\bel{q2-2}
v=e^{-mu}.
\ee
Then $v>0$ satisfies
\bel{q2-3}
-\Gd v+mv^{1-\frac 1m}=0.
\ee
There always exists a separable radial solution
\bel{q2-4}
V_m(r)=C(m,N)r^{2m}\quad\text{with }\, C(m,N)=\frac{1}{\left(2\left(2m+N-2\right)\right)^m},
\ee
hence $u(r)=\frac{1}{m}\ln \frac{C_{m,N}}{r^2}$.\smallskip

\nind {\bf1}- If $m>1$ this is an equation with a sublinear absorption. Separable non-radial solutions under the form $v(r,\gs)=r^{2m}\gw(\gs)$ ($r\geq 0$,$\gs\in S^{N-1}$) exists provided $\gw$ satisfies
\bel{q2-5}
-\Gd_g\gw-2m(2m+N-2)\gw+m\gw^{1-\frac 1m}=0\quad\text{on } S^{N-1}.
\ee
It is proved in \cite[Theorem 3.1]{BV-Gr} that if $2(2m+N-2)\leq N-1$ there exists no non-constant solution besides $C(m,N)$. Non constant separable solutions can be obtained by bifurcation methods when $2(2m+N-2)> N-1$ in \cite[Theorem 3.2]{BV-Gr}.  A description of isolated singularities is also provided in this article. The behaviour of solutions in an exterior domain can be obtained by Kelvin transform. In the case $N=2$ a complete description of separable solutions is provided in \cite {BV-Ga-Gr-V}.\smallskip

\nind {\bf 2}- If $m=1$ the function $v=e^{-u}$ satisfies
\bel{q2-6}
-\Gd v+1=0.
\ee
The radial solutions are the functions $r\mapsto \frac {r^2}{2N}+Ar^{2-N}+B$. If $u\geq 0$, then $0< v\leq 1$, hence there exists no solution in $\BBR^N$ ($u$ may change sign, then 
$A\geq 0$, note that if $A\neq 0$ these solutions are singular at $0$). Non constant separable solutions satisfy $v(r,\gs)$ are of the form $r^2\gw$ where
\bel{q2-7}
-\Gd_g\gw-2N\gw+1=0\quad\text{on } S^{N-1}.
\ee
If we set $\gf=\gw-1$, then $-\Gd_g\gf=2N\gf$. Since $2N$ is the second nonzero eigenvalue of $-\Gd_g$ on $S^{N-1}$ and $\gf_2$ an associated eigenfunction, then $\gf=c\gf_2$ and $\gw=1+c_2\gf_2$. Since $\gf_2$ is signed no separable solution $v$ satisfies $0<v\leq 1$, this means that $u$ does not stay positive.. The signed separable solutions $u$ of $(\ref{q2-6})$ are under the form 
$u=-\ln v$ where 
$$ v(r,\gs)=r^2(1+c_2\gf_2(\gs))
$$
with $|c_2|\norm{\gf_2}_{L^\infty}\leq 1$.\smallskip

\nind {\bf 3}- If $0<m<1$ we set $1-\frac 1m=-Q$, hence $Q>0$ and  $v=e^{-mu}$ satisfies
\bel{q2-8}
-\Gd v+\frac m{v^Q}=0.
\ee
The function $V_m$ defined at $(\ref{q2-4})$ is still a radial solution of the above equation. A non-radial separable solution is under the form $v(r,\gs)=r^{2m}\gw(\gs)$ where
$\gw$ satisfies 
\bel{q2-9}
-\Gd_g\gw-\tfrac{2(N(Q+1)-2Q)}{(Q+1)^2}\gw+\tfrac{1}{Q+1}\gw^{-Q}=0\quad\text{on } S^{N-1}
\ee
with $Q=\frac 1m-1\in (0,\infty)$. Constant positive solutions $\gw_s$ always exist 
\bel{q2-10}\gw_s=\left(\frac{Q+1}{2(N(Q+1)-2Q}\right)^{\frac1{Q+1}}.
\ee
If we linearise the equation at this value $\gw_s$ we obtain 
\bel{q2-10}
\CL[\gf]:=-\Gd_g\gf-\tfrac{2(N(Q+1)-2Q)}{Q+1}\gf.
\ee
A bifurcation occurs if 
\bel{q2-10}
\Gl(Q):=\tfrac{2(N(Q+1)-2Q)}{Q+1}=\gl_k=k(N+k-2)
\ee
where $\gl_k$ is a positive eigenvalue of $-\Gd_g$ on $S^{N-1}$. If $N=2$ then $\gl_1=1$ and the condition of bifurcation from $\gl_1$ is $Q=3$, that is $m=\frac{1}{4}$. If 
$N\geq 3$, then $\gl_1=N-1$. The condition of bifurcation from $\gl_1$ is $N(Q-1)=3Q-1$ which is never achieved. The condition from bifurcation from $\gl_2=2N$ is  $4Q=0$, which is also never achieved. 
Now, for $k\geq 3$, the condition of bifurcation $(\ref{q2-10})$ implies
\bel{q2-11}
-4Q=(Q+1)(k^2+(N-2)k-2N),
\ee
which cannot hold if $N\geq 3$ and any $k\geq 1$. Hence there exists  no bifurcation if $N\geq 3$ {\it We conjecture that all the positive solutions of $(\ref{q2-9})$ are constant if $N\geq 3$}.\\
Two papers are concerned with the local study of solutions of $(\ref{q2-8})$. In \cite{Ma-W} the authors study the local and global behaviour of positive solutions. Using bernstein method they provide an upper estimate of $|\nabla u|$. They also show that under suitable conditions on $N$ and $Q$ there exists no stable positive solutions and obtain radiality results under some integral constraints on the solutions.\\
 In \cite{Du-G} the authors study the equation
\bel{q2-13}
-\Gd v+\frac {|x|^\ga}{v^Q}=0.
\ee
They obtain conditions on $N$, $\ga$ and $Q>0$ for the existence or the non-existence of stable solutions. They classify the positive solutions singular at $0$.
\mysection{A priori estimate}

\subsection{Proof of Theorems 1.1 and  1.4}

For the proof of \rth{Th2} we recall the following variant of the Keller-Osserman estimate \cite[Lemma 3.1]{BV-V2}
\blemma {upp} Let $q > 1$, $d \geq 0$ and $P$ and $Q$ two continuous functions defined in $B_\gr(x)$ such
that $\inf\{P(y) : y \in B_\gr(x)\} > 0$ and $\sup\{Q(y) : y \in B_\gr(x)\} < \infty$. If $X$ is a positive $C^1$
function defined in $B_\gr(x)$ and such that
\bel{upp1}
-\Gd X+P(y)X^q\leq Q(y)+d\frac{|\nabla X|^2}{X}\quad\text{in }B_\gr(x),
\ee
then there exist positive constants $C_j= C_j(N,q,d) > 0$, $j=1,2)$ such that
\bel{upp2}
|\nabla X(x)|\leq C_1\left(\frac{1}{\gr^2\inf_{B_\gr(x)}P}\right)^{\frac 1{q-1}}+C_2\left(\sup_{B_\gr(x)}\frac{Q}{P}\right)^{\frac 1{q}}.
\ee
\es

\nind{\it Proof of \rth{Th2}.} Set $z=|\nabla u|^2$, then 
$$\BA{lll}\dsps
-\frac12\Gd z+\frac{(\Gd u)^2}{N}+\langle\nabla \Gd u,\nabla u\rangle\leq 0.
\EA$$
Hence, for $\ge>0$ small enough, 
$$\BA{lll}\dsps
-\frac12\Gd z+\frac{(mz^{\frac q2}-e^u)^2}{N}\leq e^uz +\frac {mq}2z^{\frac q2-1}\langle\nabla z,\nabla u\rangle\\[2mm]
\phantom{\dsps
-\frac12\Gd z+\frac{(mz^{\frac q2}-e^u)^2}{N}}\dsps\leq 
e^uz+\frac {mq}2\frac{|\nabla z|}{\sqrt z}z^{\frac q2}\\[2mm]
\phantom{\dsps
-\frac12\Gd z+\frac{(mz^{\frac q2}-e^u)^2}{N}}\dsps
\leq e^uz+\ge z^q+C_{\ge,m}\frac{|\nabla z|^2}{z}.
\EA$$
Then
$$\BA{lll}\dsps
-\frac12\Gd z+\frac{m^2z^{q}}{N}\leq \frac {e^{2u}}{N}+e^uz+\ge z^q+C_{\ge,m}\frac{|\nabla z|^2}{z}\\[2mm]
\phantom{\dsps
-\frac12\Gd z+\frac{m^2z^{q}}{N}}\dsps
\leq \frac {e^{2u}}{N}+\frac{m^2z^{q}}{4N}+\frac{1}{q'}\left(\frac{4N}{qm^2}\right)^{\frac{1}{q-1}}e^{q' u}+\ge z^q+C_{\ge,m}\frac{|\nabla z|^2}{z}.
\EA$$
Choosing $\ge=\frac{m^2}{4N}$, we obtain
\bel{Na-2-2}
-\frac12\Gd z+\frac{m^2z^{q}}{2N}\leq \frac {e^{2u}}{N}+\frac{1}{q'}\left(\frac{4N}{qm^2}\right)^{\frac{1}{q-1}}e^{q' u}+C_{\ge,m}\frac{|\nabla z|^2}{z}.
\ee
By applying \rlemma {upp} with $P(x)=\frac{m^2}{2N}$ and $Q(x)=\frac{e^{2u}}{N}+\frac{1}{q'}\left(\frac{4N}{qm^2}\right)^{\frac 1{q-1}}$, it follows that 
\bel{Na-2-2*}|\nabla u(x)|^2:=z(x)\leq C_1|x|^{-\frac{2}{q-1}}+C_2\max_{z\in B_{\frac{|x|}{2}}(x)}e^{\frac {2u(z)}q}+C_3\max_{z\in B_{\frac{|x|}{2}}(x)}e^{\frac {u(z)}{q-1}},
\ee
 which implies $(\ref{Na-2-1})$.
\qeda\medskip

\rcor{Th1} and \rth{Th2} follow immediately from $(\ref{Na-2-1})$. The proof of \rth{Th3} is much more elaborate.\medskip

\nind{\it Proof of \rth{Th3}}. The estimates $(\ref{Na-20-4})$ are equivalent to 
\bel{Na-2-3}
e^{\frac{u(x)}{q}}+|\nabla u(x)|\leq \frac C{|x|}\quad\text{for all }x\in B_{\frac \gr2}\setminus\{0\}.
\ee
Set $w=e^u$, then $w$ is positive and satisfies
\bel{Na-2-4}
-\Gd w+\frac{|\nabla w|^2}{w}+m\frac{|\nabla w|^q}{w^{q-1}}-w^2=0.
\ee
We will prove $(\ref{Na-2-3})$ under the equivalent form
\bel{Na-2-5}
w^\frac{1}{q}(x)+\frac{|\nabla w(x)|}{w(x)}\leq \frac C{|x|}\quad\text{for all }x\in B_{\frac \gr2}\setminus\{0\}.
\ee
We set 
$$M(x):=e^{\frac{u(x)}{q}}=w^\frac{1}{q}(x),
$$
and we claim that $(\ref{Na-2-5})$ will hold if 
\bel{Na-2-5'}
M(x)\leq \frac C{|x|}\quad\text{for all }x\in B_{\frac \gr 2}\setminus\{0\}.
\ee
Indeed, if this estimate holds, we have from $(\ref{Na-2-1})$,
\bel{Na-2-5''}
\max_{B_{\frac {|x|}4}(x)}|\nabla u(x)|\leq C\left(|x|^{-\frac{1}{q-1}}+|x|^{-1}+|x|^{-\frac{q}{2(q-1)}}\right)\leq 2C|x|^{-1}\quad\text{for all }x\in B_{\frac \gr2}\setminus\{0\}.
\ee
We apply the doubling Lemma \cite{PQS} with $X=\overline {B_1}$, $D=\overline {B_{\frac 12}}\setminus\{0\}$, $\Gs=\overline {B_{\frac 12}}$ and $\Gg=\{0\}$. Let
$k > 0$: If $M$ is any function defined in $D$ bounded on compact subsets of $D$; and $y\in D$ such that
$M(y)|y|>2k$, then there exists $x$ such that $0 < |x|\leq\frac 12$ and
\bel{Na-2-5'''}\left\{\BA{lll}
M(x)|x|>2k\\[1mm]
M(x)\geq M(y)\\
\dsps M(z)\leq 2M(x)\quad\text{for all }z\text{ s.t. }|z-x|\leq\frac{k}{M(x)}.
\EA\right.
\ee
Assume that $(\ref{Na-2-5'})$ does not hold. Then there exists $y_n\in \overline {B_{\frac 12}}\setminus\{0\}$ such that $M(y_n)|y_n|>2n$. Then there exists $x_n\in \overline {B_{\frac 12}}\setminus\{0\}$ such that 
\bel{Na-2*-5'''}\left\{\BA{lll}
M(x_n)|x_n|>2n\\[1mm]
\phantom{|x_n|}M(x_n)\geq M(y_n)\\
\phantom{|x_N|}\dsps M(Z)\leq 2M(x_n)\quad\text{for all }Z\text{ s.t. }|Z-x|\leq\frac{n}{M(x_n)}.
\EA\right.
\ee
We can extract a subsequence such that $\dsps\lim_{n\to\infty}x_n = 0$ (indeed if $\dsps\liminf_{n\to\infty}|x_n| > 0$,
since $M(x_n)$ is bounded on compact set of  $\overline {B_{\frac 12}}\setminus\{0\}$, this is in contradiction with $M(x_n)|x_n|>2n$).  Thus $\lim_{n\to\infty}M(x_n)=\infty$. Next we  set 
$$w_n(x)=\frac{w(Z(x,n))}{M^q(x_n)}\quad\text{where }Z(x,n)=x_n+\frac{x}{M(x_n)}.
$$
Then $w_n(0)=1$ and $w_n\leq 2$ on $B_n$. But we have
$$\nabla w_n(x)=\frac{\nabla w(Z(x,n))}{M^{q+1}(x_n)}\,,\; \Gd w_n(x)=\frac{\Gd w(Z(x,n))}{M^{q+2}(x_n)}\,\text{ and }\; \frac{|\nabla w_n(x)|}{w_n(x)}
=\frac{1}{M(x_n)}\frac{|\nabla w(Z(x,n))|}{w(Z(x,n))}.
$$
Thus form $(\ref{Na-2-4})$
$$\BA{lll}\dsps
-\Gd w_n(x)=\frac{1}{M^{q+2}(x_n)}\left(w^2-\frac{|\nabla w|^2}{w}-m\frac{|\nabla w|^q}{w^{q-1}}\right)(Z(x,n))\\[4mm]
\phantom{-\Gd w_n(x)}\dsps=\frac{1}{M^{q+2}(x_n)}\left(M^{2q}(x_n)w^2_n(x)-M^{q+2}(x_n)\frac{|\nabla w|^2}{w}-M^{2q}(x_n)m\frac{|\nabla w|^q}{w^{q-1}}\right)(Z(x,n))\\[4mm]
\phantom{-\Gd w_n(x)}\dsps=M^{q-2}(x_n)\left(w_n^2(x)-M^{2-q}(x_n)\frac{|\nabla w_n(x)|^2}{w_n(x)}-m\frac{|\nabla w_n(x)|^q}{w^{q-1}_n(x)}\right).
\EA
$$
But $M^{2-q}(x_n)=\ge_n\to 0$ since $q>2$, and
$$M^{2-q}(x_n)\left(-\Gd w_n(x)+\frac{|\nabla w_n(x)|^2}{w_n(x)}\right)=w^2_n(x)-m\frac{|\nabla w_n(x)|^q}{w^{q-1}_n(x)}.
$$
Thus 
$$M^{2-q}(x_n)w^{q-1}_n(x)\left(-\Gd w_n(x)+\frac{|\nabla w_n(x)|^2}{w_n(x)}\right)=w^{q+1}_n-m|\nabla w_n(x)|^q.
$$
Put $a_n(x)=M^{2-q}(x_n)w^{q-1}_n(x)$. Then $\dsps\sup_{B_n}a_n(x)\to 0$ when $n\to\infty$ and we have
equivalently
\bel{Na-2-6}
-a_n(x)\Gd w_n(x)+M^{2-q}(x_n)w^{q-2}_n(x)|\nabla w_n(x)|^2=w_n^{q+1}(x)-m|\nabla w_n(x)|^q.
\ee
For the estimates of the gradient, we know from $(\ref{Na-2-2*})$ that for a solution defined in $B_\gr(a)$ and since $q>2$  that
$$\sup_{B_{\frac \gr 2}(a)}\frac{\nabla w}{w}\leq c_1\gr^{-\frac 1{q-1}}+c'_2\sup_{B_{\gr}(a)}w^{\frac 1q}<\infty.
$$
Hence with $a=Z(x,n)=x_n+\frac x{M(x_n)}$ and $\gr=\frac{n}{M(x_n)}$,
$$\sup_{B_{\frac{n}{2M(x_n)}}(x_n)}\frac{\nabla w(Z)}{w(Z)}\leq c_1\left(\frac{2M(x_n)}{n}\right)^{\frac1{q-1}}+c'_2\sup_{B_{\frac{n}{M(x_n)}}(x_n)}w^{\frac 1q}.
$$
Equivalently
$$M(x_n)\sup_{B_{\frac n2}}\frac{|\nabla w_n|}{w_n}\leq c_1\left(\frac{2M(x_n)}{n}\right)^{\frac1{q-1}}M^{\frac{2-q}{q-1}}(x_n)+c'_2M(x_n)\sup_{B_{\frac{n}{M(x_n)}}(x_n)}w_n^{\frac 1q},
$$
and finally
$$\sup_{B_{\frac n2}}\frac{|\nabla w_n|}{w_n}\leq c'_1M^{\frac{2-q}{q-1}}(x_n)+c'_2\sup_{B_{\frac{n}{M(x_n)}}(x_n)}w_n^{\frac 1q},
$$
But this last quantity is upper bounded since $M^{2-q}(x_n)\to 0$ and $w_n\leq 2$ in $B_n$. Therefore $\frac{|\nabla w_n|}{w_n}$ is bounded 
in $B_{\frac n2}$ and thus it is the same with $|\nabla w_n|$. Therefore we write $(\ref{Na-2-6})$ as
\bel{Na-2-7}
-a_n(x)\Gd w_n(x)=B_n(x)+w_n^{q+1}(x)-m|\nabla w_n(x)|^q,
\ee
with $B_n(x)=M^{2-q}(x_n)w_n^{q-2}(x)|\nabla w_n(x)|^2\to 0$ locally uniformly. Then there exists a subsequence still denoted $w_n$ and a function $W\in W^{1,\infty}(\BBR^N)$ such that $w_n\to W$ locally uniformly in $\BBR^N$ and $\nabla w_n\to\nabla W$ in the weak-star topology of $L_{loc}^{\infty}(\BBR^N)$. Even if $a_n$ is not a constant, it converges to $0$ uniformly in $B_n$. By a standard adaptation of 
\cite[Proposition IV.1]{CrLi}, we conclude that W is a viscosity solution of the eikonal equation
$$W^{q+1}-m|\nabla W|^q=0
$$
in $\BBR^N$. By \cite[Proposition IV.3]{CrLi} the only bounded solution is $W=0$. Because $w_n(0)=1$ it follows that  $W(0)=1$, a contradiction.\qeda

\subsection{Estimates on supersolutions}

We first give a monotonicity result which is a corrected version of \cite[Lemma 2.2]{BV-V2}.

\blemma{erra}Let $N\geq 2$ and $q>1$. If $u\in C^2(B_{\gr_0}\setminus\{0\})$ is a positive function such that 
\bel{Er1}-\Gd u+|\nabla u|^q\geq 0\ee
in $B_{\gr_0}\setminus\{0\}$. Then there exists $\gr_1\in (0,\gr]$ such that the function $\gm(r)=\min_{|x|=r}u(x)$ is monotone on $(0,\gr_1)$. Furthermore, if $1<q\leq 2$ we have that $\gr_1=\gr_0$ and 
$\gm$ is nonincreasing on $(0,\gr_0]$.
If  $u\in C^2(B^c_{\gr_0})$ is a positive function satisfying the above inequality in $B^c_{\gr_0}$,
there exists  $\gr_1\geq \gr_0$ such that $\gm$  is monotone on $[\gr_1,\infty)$.
\es

\nind\Proof  The case of an exterior domain is treated in \cite[Lemma 5]{Al-Me-Qu}.  In the case of the punctured ball, for any $0<s<t<\gr$ we have by the comparison principle  
$$\dsps \min_{ s\leq|x|\leq t} u(x)=\min_{_{ s\leq r\leq t}}\gm(r)=\min\{\gm(s),\gm(t)\},$$ 
thus the function $\dsps(s,t)\mapsto \min\{\gm(s),\gm(t)\}$ is nondecreasing in $s$ and nonincreasing in $t$, hence $\gm$ cannot have any local strict minimum.\\
If the function  $\gm$ which is upper semicontinuous,  is bounded on $(0,\gr]$ its supremum is achieved at $\gr_1\in [0,\gr]$. If $0<\gr_1<\gr$, then 
$\gm(s)=\min\{\gm(s),\gm(r_1)\}$ and it is nondecreasing in $s$. If $\gr_1=0$ (resp. $\gr_1=\gr$ ), then by monotonicity $\gm(s)=\min\{\gm(0_+),\gm(s)\}$ is nonincreasing in $s$ on $(0,\gr]$ (resp. $\gm(s)=\min\{\gm(s),\gm(r_0)\}$ is nondecreasing in $s$).\\
If $\gm$ is unbounded, there exists $\{r_n\}$ decreasing to $0$ such that $\gm(r_n)\to\infty$ when $n\to\infty$. Hence for $0<r<s<r_0$ there exists $r_n$ such that $r_n<r$ and 
$\gm(r_n)\geq max\{\gm(r),\gm(s)\}$. Then $\gm(r)=\min\{\gm(r_n),\gm(r)\}\geq \min\{\gm(r_n),\gm(s)\}=\gm(s)$.\\
Finally, in the case $1<q\leq 2$ we proceed by contradiction in assuming that there exist $\gr_1\in (0,\gr)$ such that $\gm$ is nondecreasing on $(0,\gr_1)$ and $\gm(0)=\gth\geq 0$. The function $h_c(r)=cr$ is a radial subsolution of the equation in $B_{\gr_1}\setminus\{0\}$ smaller that $u$ on $\prt B_{\gr_1}$ and near $0$ (up to changing $u$ in $u+\ge$) provided $0<c\leq \min\left\{\frac{\gm(\gr_1)}{\gr_1},\left(\frac{N-1}{\gr_1}\right)^{\frac1{q-1}}\right\}$. Hence there exists 	a positive radial function $h$ satisfying $0\leq h(r)\leq \gm(r) $ and 
\bel{Er2}-\Gd h+|\nabla h|^q=0.\ee
By an explicit integration it is easy to verify that when $1<q\leq 2$ there exists no positive bounded radial solution (even supersolution) $h$ in $B_{\gr_1}\setminus\{0\})$ of $(\ref{Er2})$ which is  a contradiction.
$\phantom{-----------}$
\qeda\medskip

\nind {\it Proof of \rth{Th4}}. We set $w=e^u$, hence 
\bel{Na-2-8}
-\Gd w+\frac{|\nabla w|^2}{w}+m\frac{|\nabla w|^q}{w^{q-1}}\geq w^2.
\ee
Set
$$\gm(r)=\min_{|x|=r}u(x)\quad\text{and }\;M(r)=\min_{|x|=r}w(x)=e^{\gm(r)}.
$$
By \rlemma{erra} applied to $u$ there exists $\gr_1\in [0,\gr_0]$ such that the function $\gm$ is monotone on $(0,\gr_1]$, and thus it is the same with $M$. 
If $M$ is nondecreasing on $(0,\gr_1]$, estimate $(\ref{Na-22})$ holds true with $C=\norm{w}_{L^{\infty}(B_{\gr_0})}$.   Assume now that $\gm$ is nonincreasing on $(0,\gr_1]$. Let $0<\ge\leq \frac 12$, $0<\gr<\frac {\gr_1}2$ and
 $\phi_\ge$ be a radial cut-off function such that 
$\phi_\ge=0$ on $[0,1-\ge]\cup [1+\ge,\infty)$, $\phi_\ge=1$ on  $[1-\frac\ge2,1+\frac\ge2]$ and $|\phi_\ge'|\leq \frac C\ge$ on $[1-\ge,1-\frac \ge2]\cup [1+\frac \ge2,1+\ge]$. We consider the function
$$v(x)=w(x)-M(\gr)\phi_\ge\left(\frac{|x|}{\gr}\right).
$$
There exists a point $x_{\gr,\ge}$ such that $|x_{\gr,\ge}|=\gr$ where $w(x_{\gr,\ge})=M(\gr)$ and $|x_{\gr,\ge}|=\gr$, thus $v(x_{\gr,\ge})=0$; since $w$ is defined in $B_{\gr_0}\setminus\{0\}$, we have 
$v=w$ in $\left(B_{(1-\ge)\gr}\setminus\{0\}\right)\cup \left(B_{\gr_1}\setminus \overline B_{(1+\ge)\gr}\right)$. Therefore $v$ achieves its
nonpositive minimum at some point $\tilde x_{\gr,\ge}\in B_{(1+\ge)\gr}\setminus \overline B_{(1-\ge)\gr}$ where $\nabla v(\tilde x_{\gr,\ge})=0$ and $\Gd v(\tilde x_{\gr,\ge})\geq 0$ with $v(\tilde x_{\gr,\ge})\leq 0$. Then  $w(\tilde x_{\gr,\ge})\leq M(\gr)$, so that $M(|\tilde x_{\gr,\ge}|)\leq M(\gr)$. There holds
$\nabla w(\tilde x_{\gr,\ge})=M(\gr)\nabla \phi_\ge\left(\frac {|\tilde x_{\gr,\ge}|}\gr\right)$ and $\Gd w(\tilde x_{\gr,\ge})\geq M(\gr)\Gd \phi_\ge\left(\frac {|\tilde x_{\gr,\ge}|}\gr\right)$. Then there exists $C=C(N,q,m)>0$ such that 
\bel{Na-2-9}\BA{lll}\dsps
w^2(\tilde x_{\gr,\ge})\leq -\Gd w(\tilde x_{\gr,\ge})+\frac{|\nabla w(\tilde x_{\gr,\ge})|^2}{w(\tilde x_{\gr,\ge})}+m\frac{|\nabla w(\tilde x_{\gr,\ge})|^q}{w^{q-1}(\tilde x_{\gr,\ge})}
\\[4mm]
\phantom{w^2(\tilde x_{\gr,\ge})}\dsps \leq -M(\gr)\Gd\phi_\ge\left(\frac {|\tilde x_{\gr,\ge}|}\gr\right)+
M^2(\gr)\frac{\left|\nabla\phi_\ge\left(\frac{|\tilde x_{\gr,\ge}|}{\gr}\right)\right|^2}{w(\tilde x_{\gr,\ge})}+mM^q(\gr)\frac{\left|\nabla\phi_\ge\left(\frac{|\tilde x_{\gr,\ge}|}{\gr}\right)\right|^q}{w^{q-1}(\tilde x_{\gr,\ge})}\\[4mm]
\phantom{w^2(\tilde x_{\gr,\ge})}\dsps 
\leq C\left(\frac{M(\gr)}{\ge^2\gr^2}+\frac{M^2(\gr)}{\ge^2\gr^2w(\tilde x_{\gr,\ge})}+\frac{ M^q(\gr)}{\ge^q\gr^qw^{q-1}(\tilde x_{\gr,\ge})}\right).
\EA\ee
If $q>2$ we obtain with $C=C(N,q,m)>0$ varying from one occurence to another, 
$$\BA{lll}\dsps
w^{q+1}(\tilde x_{\gr,\ge})\leq  C\left(\frac{M(\gr)}{\ge^2\gr^2}w^{q-1}(\tilde x_{\gr,\ge})+\frac{M^2(\gr)}{\ge^2\gr^2}w^{q-2}(\tilde x_{\gr,\ge})+\frac{ M^q(\gr)}{\ge^q\gr^q}\right)\\[4mm]
\phantom{w^{q+1}(\tilde x_{\gr,\ge})}\dsps 
\leq C\left(\frac{M(\gr)}{\ge^2\gr^2}M^{q-1}(\gr)+\frac{M^2(\gr)}{\ge^2\gr^2}M^{q-2}(\gr)+\frac{ M^q(\gr)}{\ge^q\gr^q}\right)\\[4mm]
\phantom{w^{q+1}(\tilde x_{\gr,\ge})}\dsps
\leq C\frac{M^q(\gr)}{\ge^q\gr^q},
\EA$$
and $\dsps w^{q+1}(\tilde x_{\gr,\ge})\geq \min_{(1-\ge)\gr\leq r\leq (1+\ge)\gr}M^{q+1}(r)\geq M^{q+1}((1+\ge)\gr)$. Then we get the estimate
\bel{Na-2-10}
M((1+\ge)\gr)\leq C\frac{M^{\frac{q}{q+1}}(\gr)}{\ge^{\frac{q}{q+1}}\gr^{\frac{q}{q+1}}}.
\ee
Applying the bootstrap argument of \cite[Lemma 2.1]{BV-V2} with $\Gf(t)=t^{-\frac{q}{q+1}}$ and $d=\frac q{q+1}$, we obtain 
\bel{Na-2-11}
M(\gr)\leq C\left(\Gf(\gr)\right)^{\frac{1}{1-d}}=C\gr^{-q}, 
\ee
which is the claim in that case. \\
If $q<2$ we obtain from $(\ref{Na-2-9})$, 
$$\BA{lll}\dsps
w^{3}(\tilde x_{\gr,\ge})\leq  C\left(\frac{M(\gr)}{\ge^2\gr^2}w(\tilde x_{\gr,\ge})+\frac{M^2(\gr)}{\ge^2\gr^2}+\frac{ M^q(\gr)}{\ge^q\gr^q}w^{2-q}(\tilde x_{\gr,\ge})\right)\\[4mm]
\phantom{w^{3}(\tilde x_{\gr,\ge})}\dsps 
\leq C\left(\frac{M(\gr)}{\ge^2\gr^2}M(\gr)+\frac{M^2(\gr)}{\ge^2\gr^2}+\frac{ M^q(\gr)}{\ge^q\gr^q}M^{2-q}(\gr)\right)\\[4mm]
\phantom{w^{3}(\tilde x_{\gr,\ge})}\dsps
\leq C\frac{M^2(\gr)}{\ge^2\gr^2},
\EA$$
and $w^3(\tilde x_{\gr,\ge})\geq M^3((1+\ge)\gr)$. We apply again \cite[Lemma 2.1]{BV-V2} with $\Gf(t)=t^{-\frac{2}{3}}$ and $d=\frac 2{3}$ and deduce that 
\bel{Na-2-12}
M(\gr)\leq C\left(\Gf(\gr)\right)^{\frac{1}{1-d}}=C\gr^{-2}, 
\ee
which ends the proof.
\qeda\medskip
 
Following the ideas introduced in \cite{BiPo} and adapted to some equations with gradient term by \cite{Rui} we prove the following,


\bprop{Lemm1} Let $N>2>q>1$, $\ell>\max\{\frac q{2-q},\frac N{N-2}\}$ and $\gn>0$. Assume that $u\in C^1(\overline B_{\gr_0}\setminus\{0\})$ is a function which Laplacian in the sense of distributions  belongs to 
$L^1_{loc}(B_{\gr_0}\setminus\{0\})$ and such that
$u\geq -a$ for some real number $a$. If there holds 
\bel{E4-4}\dsps
-\Gd u+m|\nabla u|^q\geq F:=A(u+a+\gn)^\ell
\ee
 in $\CD'(B_{\gr_0}\setminus\{0\})$, then $u\in L^s(B_{\gr_0})$ for all $1<s<\ell$ and $|\nabla u|\in L^\gs(B_{\gr_0})$ for all $1<\gs<\frac{2\ell}{\ell-1}$. 
\es
\Proof Notice that the assumption on $\ell$ is equivalent to $q<\frac{2\ell}{\ell+1}$ (subcriticality assumption). Without loss of generality we can assume that $a=0$. For $\gr\in (0,\frac {\gr_0}2)$, let $\xi_\gr\in C^{\infty}_c(B_{\gr_0}\setminus\{0\})$ with values in $[0,1]$ such that $\xi_\gr=0$ in $B_{\frac\gr 4}\cup B^c_{2\gr}$, $\xi_\gr=1$ in $B_{\gr}\setminus B_{\frac\gr 2}$ and 
$|\nabla \xi_\gr|\leq C\gr^{-1}{\bf I}_{B_{2\gr}\setminus B_{\frac\gr 4}}$. We take $\phi_\gr=\xi_\gr^\gl$  for test function where $\gl>0$ is large enough and will be made precise later on. If $-1<\ga<0$ we multiply $(\ref{E4-4})$ by $(u+\gn)^{\ga}\xi^\gl_\gr$. We set $u_\gn=u+\gn$ and obtain after integration
$$\BA{lll}\dsps A\int_{B_{\gr_0}}u_\gn^{\ell+\ga}\xi_\gr^\gl dx+|\ga|\int_{B_{\gr_0}}u_\gn^{\ga-1}\xi_\gr^\gl |\nabla u_\gn|^2dx\\[4mm]
\phantom{--------}\dsps\leq
 \int_{B_{\gr_0}}u_\gn^\ga\xi_\gr^{\gl}|\nabla u_\gn|^qdx+\gl \int_{B_{\gr_0}}u_\gn^\ga\xi_\gr^{\gl-1}\nabla u_\gn.\nabla\xi_\gr dx.
\EA$$
For any $\ge>0$ small enough, we have by Cauchy-Schwarz inequality,
$$\left|\int_{B_{\gr_0}}u_\gn^\ga\xi_\gr^{\gl-1}\nabla u_\gn.\nabla\xi_\gr dx\right|\leq \ge\int_{B_{\gr_0}}u_\gn^{\ga-1}\xi_\gr^\gl |\nabla u_\gn|^2dx+\frac{1}{\ge}
\int_{B_{\gr_0}}u_\gn^{\ga+1}\xi_\gr^{\gl-2}|\nabla \xi_\gr|^2 dx,
$$
and
$$\BA{lll}\dsps 
 \int_{B_{\gr_0}}u_\gn^\ga\xi_\gr^{\gl}|\nabla u_\gn|^qdx\leq   \ge\int_{B_{\gr_0}}u_\gn^{\ga-1}\xi_\gr^\gl |\nabla u_\gn|^2dx+\frac {1}{\ge}\int_{B_{\gr_0}}u_\gn^{\ga+\frac q{2-q}}\xi_\gr^\gl  dx\\[4mm]
 \phantom{\dsps 
 \int_{B_{\gr_0}}u_\gn^\ga\xi_\gr^{\gl}|\nabla u^+_\gn|^qdx}
 \dsps \leq 
 \ge\int_{B_{\gr_0}}u_\gn^{\ga-1}\xi_\gr^\gl |\nabla u_\gn|^2dx+\ge\int_{B_{\gr_0}}u_\gn^{\ga+\ell}\xi_\gr^{\gl}+C_\ge\gr^N,
\EA$$
from the choice of $\xi_\gr$. Choosing $\ge$ small enough we obtain, with $\gk>1$ and $\gk'=\frac{\gk}{\gk-1}$,
$$\frac A2\int_{B_{\gr_0}}u_\gn^{\ell+\ga}\xi_\gr^\gl dx+\frac{|\ga|}2\int_{B_{\gr_0}}u_\gn^{\ga-1}\xi_\gr^\gl |\nabla u_\gn|^2dx\leq C'_\ge\gr^N+C'_\ge\int_{B_{\gr_0}}u_\gn^{\ga+1}\xi_\gr^{\gl-2}|\nabla \xi_\gr|^2 dx,
$$
$$\BA{lll}\dsps\int_{B_{\gr_0}}u_\gn^{\ga+1}\xi_\gr^{\gl-2}|\nabla \xi_\gr|^2 dx=\int_{B_{\gr_0}}\xi_\gr^{\gth}u_\gn^{\ga+1}\xi_\gr^{\gl-2-\gth}|\nabla \xi_\gr|^2 dx\\[4mm]
\phantom{\dsps\int_{B_{\gr_0}}u_\gn^{\ga+1}\xi_\gr^{\gl-2}|\nabla \xi_\gr|^2 dx}\dsps\leq \left(\int_{B_{\gr_0}}\xi_\gr^{\gk\gth}u_\gn^{\gk(\ga+1)}dx\right)^{\frac1\gk}
\left(\int_{B_{\gr_0}}\xi_\gr^{(\gl-2-\gth)\gk'}|\nabla \xi_\gr|^{2\gk'} dx\right)^{\frac1\gk'},
\EA$$
then
\bel{E4-5}\BA{lll}\dsps\frac A2\int_{B_{\gr_0}}u_\gn^{\ell+\ga}\xi_\gr^\gl dx+\frac{|\ga|}2\int_{B_{\gr_0}}u_\gn^{\ga-1}\xi_\gr^\gl |\nabla u_\gn|^2dx\\[4mm]
\phantom{------}\dsps\leq C'_\ge\gr^N+C'_\ge\left(\int_{B_{\gr_0}}\xi_\gr^{\gk\gth}u_\gn^{\gk(\ga+1)}dx\right)^{\frac1\gk}
\left(\int_{B_{\gr_0}}\xi_\gr^{(\gl-2-\gth)\gk'}|\nabla \xi_\gr|^{2\gk'} dx\right)^{\frac1{\gk'}}.
\EA\ee
Fix $\gk=\frac{\ell+\ga}{1+\ga}$, then $\gk'=\frac{\ell+\ga}{\ell-1}$ and thus $\gk(\ga+1)=\ell+\ga$. Take $\gk\gth=\gl-2$, hence 
$\gl-2-\gth=(\gk-1)\gth$, hence $(\gl-2-\gth)\gk'=\gth\gk$.
Therefore
\bel{E4-6}
\frac A2\int_{B_{\gr_0}}u_\gn^{\ell+\ga}\xi_\gr^\gl dx+\frac{|\ga|}2\int_{B_{\gr_0}}u_\gn^{\ga-1}\xi_\gr^\gl |\nabla u_\gn|^2dx\leq 
C'_\ge\gr^N+C^*_\ge\gr^{\frac{N(\ell-1)}{\ell+\ga}-2}\left(\int_{B_{\gr_0}}u_\gn^{\ell+\ga}\xi_\gr^\gl dx\right)^{\frac 1\gk}.
\ee
If we set $X=\left(\int_{B_{\gr_0}}u_\gn^{\ell+\ga}\xi_\gr^\gl dx\right)^{\frac 1\gk}$ and put
$$\Gth(X)=\frac A2X^\gk-C^*_\ge\gr^{\frac{N(\ell-1)}{\ell+\ga}-2}X-C'_\ge\gr^N,
$$
it is straightforward to see that
$$\Gth(X)\leq 0\Longrightarrow X^\gk\leq \tilde C_\ge\left(\gr^N+\gr^{N-2\frac{\ell+\ga}{\ell-1}}\right)\leq 2\tilde C_\ge\gr^{N-2\frac{\ell+\ga}{\ell-1}}.
$$
This yields
\bel{E4-7}
\int_{\frac\gr 2\leq|x|\leq\gr}u_\gn^{\ell+\ga} dx\leq2\tilde C_\ge\gr^{N-2\frac{\ell+\ga}{\ell-1}}.
\ee
From $(\ref{E4-6})$ we also obtain an estimate of the gradient
\bel{E4-8}\BA{lll}\dsps
\frac{|\ga|}2\int_{B_{\gr_0}}u_\gn^{\ga-1}\xi_\gr^\gl |\nabla u_\gn|^2dx\leq \hat C_\ge\gr^{N-2\frac{\ell+\ga}{\ell-1}}.
\EA\ee
We can improve this estimate in the following way. For $\gs>0$ and $b>1$ which will be specified
$$\BA{lll}\dsps\int_{B_{\gr_0}}\xi_\gr^\gl |\nabla u_\gn|^\gs dx=\int_{B_{\gr_0}}\xi_\gr^\gl u_\gn^\gd u_\gn^{-\gd} |\nabla u_\gn|^\gs dx\\[4mm]
\phantom{-----------}\dsps\leq \left(\int_{B_{\gr_0}} u_\gn^{\gd b} \xi_\gr^\gl dx\right)^\frac{1}{b}\left(\int_{B_{\gr_0}} u_\gn^{-\gd b'} |\nabla u_\gn|^{\gs b'}\xi_\gr^\gl dx\right)^\frac{1}{b'}.
\EA$$
Choosing $b=\frac{\ell+1}{1-\ga}$, hence $b'=\frac{\ell+1}{\ell+\ga}$, $\gs=2\frac{\ell+\ga}{\ell+1}$ and $\gd=\frac{(\ell+\ga)(1-\ga)}{\ell+1}$ we obtain, since $\gs b'=2$,
\bel{E4-9}\int_{\frac{\gr}{2}<|x|<\gr}|\nabla u_\gn|^{2\frac{\ell+\ga}{\ell-1}}dx\leq \left(\int_{B_{\gr_0}}u_\gn^{\ga-1}|\nabla u_\gn|^{2}\xi^\gl_\gr dx\right)^\frac{\ell+\ga}{\ell+1}\left(\int_{B_{\gr_0}}u_\gn^{\ell+\ga}\xi^\gl_\gr dx\right)^\frac{1-\ga}{\ell+1}\leq C^{\#}_\ge\gr^{N-2\frac{\ell+\ga}{\ell-1}}.
\ee
We replace now $\gr$ by $2^{-n}\gr$ ($n\in\BBN$) and get
\bel{E4-10}
\int_{|x|\leq\gr}u_\gn^{\ell+\ga} dx=\sum_{n=0}^{\infty}\int_{2^{-n-1}\gr\leq |x|\leq 2^{-n}\gr}u_\gn^{\ell+\ga} dx\leq 2\tilde C_\ge\gr^{N-2\frac{\ell+\ga}{\ell-1}}\sum_{n=0}^{\infty}2^{-n}
= 4\tilde C_\ge\gr^{N-2\frac{\ell+\ga}{\ell-1}}.
\ee
Similarly
\bel{E4-11}
\int_{|x|\leq\gr}|\nabla u_\gn|^{2\frac{\ell+\ga}{\ell-1}}dx\leq  2\tilde C^{\#}_\ge\gr^{N-2\frac{\ell+\ga}{\ell-1}}.
\ee
Since we can take $\ga$ as close to $0$ as we want, we obtain the results.
\qeda\medskip

\nind\Remark Estimates $(\ref{E4-7})$ and $(\ref{E4-8})$ are valid if $N=2$ since they hold under the mere assumption $\ell>\frac q{2-q}$.

\bcor{Cor1} Assume $N>2>q>1$. If $u\in C^1(\overline B_{\gr_0}\setminus\{0\})$, with $\Gd u\in L^1_{loc}(B_{\gr_0}\setminus\{0\})$ in the sense of distributions, is a function bounded from below, satisfying 
\bel{E4-4*}\dsps
-\Gd u+m|\nabla u|^q\geq e^{u}
\ee
 in $\CD'(B_{\gr_0}\setminus\{0\})$, then $e^u\in L^1(B_{\gr_0})$, $|\nabla u|\in L^\gs(B_{\gr_0})$ for all $1<\gs<2$ and there exists $\gg\geq 0$ and $G\in L^{1}_{loc}(B_{\gr_0})$ such that 
   \begin{equation}\label {Na-2-13}
 -\Delta u+m|\nabla u|^q= e^{u}+G+c_N \gg\gd_0\quad\text{in }\CD'(B_{\gr_0}).
 \end{equation}
\es
\Proof Since $u\geq \gn$, the function $u_\gn=u-\gn$ is nonnegative and satisfies
$$-\Gd u_\gn+m|\nabla u_\gn|^q\geq e^{\gn}e^{u_\gn}
$$
Up to changing $m$ and ${\gr_0}$ into $\tilde m>0$ and $\tilde{\gr_0}>0$ respectively by using the scaling transformation defined in  $(\ref{Na-12})$, we can assume that $\gn=0$. We apply \rprop{Lemm1} with $\ell$ as large as we want, we obtain that for any $\tilde q<2$ one has
\bel{E4-11**}
\int_{|x|\leq\gr}|\nabla u|^{\tilde q}dx\leq  2\tilde C^{\#}_\ge\gr^{N-\tilde q}.
\ee
In particular, this implies that $|\nabla u|^q\in L_{loc}^1(B_{\tilde\gr})$. By the Brezis-Lions \cite{Br-Li} result it implies that $\Gd u$ and $e^u$ are locally integrable in $B_{\tilde\gr}$ and there exists $\gg\geq 0$ and $G\in L^{1}_{loc}(B_{\gr_0})$  such that 
 $( \ref{Na-2-13})$ is verified.
\qeda\medskip


\subsection{Estimates on solutions}
We first consider solutions of $(\ref{Na-1})$ in $B_\gr\setminus\{0\}$ which are positive near $0$. The first observation with such solutions is the following 

\bprop{Dirac2} Let $N\geq 2$, $q>1$ and $u\in C^2(\overline B_{\gr_0}\setminus\{0\})$ be a positive solution of $(\ref{Na-1})$ in $B_{\gr_0}\setminus\{0\}$. Then $e^u$ is integrable in 
 $B_{\gr_0}$ if and only if  $|\nabla u|^q$ is integrable too in $B_{\gr_0}$.  Furthermore  there exists $\gg\geq 0$ such that $(\ref{Na-2-13})$ holds and we have that
 $\gg\in [0,2]$ if $N=2$ and $\gg=0$ if $N\geq 3$.
\smallskip
\es
\Proof (i)- If $|\nabla u|^q\in L^1(B_{\gr_0})$, the integrability of $e^u$ in $L^1(B_{\gr_0})$ follows by the result by \cite{Br-Li}. \\
(ii)- Let us assume now that $e^u\in L^1(B_{\gr_0})$. We proceed by contradiction  in supposing that $|\nabla u|^q\notin L^1(B_{\gr_0})$. If we denote by $(r,\gs)\in \BBR_+\ti S^{N-1}$ 
(with $\BBR_+=[0,\infty)$) the spherical coordinates in $\BBR^N$, we have, by averaging the equation on $S^{N-1}$, the following relation
$$-(r^{N-1}\bar u_r)_r=r^{N-1}\left(\overline {e^u}-m\overline{|\nabla u|^q}\right):=r^{N-1}F(r),
$$
which yields
\bel{E4-11+}-\gr_0^{1-N}\bar u_r(\gr_0)+r^{N-1}\bar u_r(r)=\int_r^{\gr_0}s^{N-1}F(s)ds.
\ee
Since $e^u\in L^1(B_{\gr_0})$ and  $|\nabla u|^q\notin L^1(B_{\gr_0})$, there holds
$$\lim_{r\to 0}\int_r^{\gr_0}s^{N-1}F(s)ds=-\infty.
$$
Therefore
$$\lim_{r\to 0}r^{N-1}\bar u_r(r)=-\infty\quad\text{which implies that }\, \lim_{r\to 0}\frac{\bar u(r)}{H_N(r)}=\infty,
$$
where we denote by $H_N$ the fundamental solution of $-\Gd$ in $\BBR^N$. This implies that 
$$\int_{B_{\gr_0}}e^{\bar u}dx=\infty.
$$
Since by convexity $\overline{e^{ u}}\geq e^{\bar u}$, it follows
$$\int_{B_{\gr_0}}e^{ u}dx=\infty,
$$
which contradicts the assumption (ii). Hence   $(\ref{Na-2-13})$ holds for some $\gg\geq 0$. \\
When $N\geq 3$ we have from $(\ref{E4-11+})$, that $r^{N-1}\bar u_r(r)$ admits a finite limit $\ell$ when $r\to 0$. If $\ell\neq 0$, then it should be negative and 
$$\bar u(r)\sim\frac{\ell}{2-N}r^{2-N}\quad\text{as }r\to 0.
$$
Now $\overline{e^u}\geq e^{\bar u}$, therefore the previous relation implies  $e^u\notin L^1(B_{\gr_0})$ which is again  in contradiction with (ii). Thus $\ell=0$ which implies that $\gg=0$.\\
When $N=2$ and $\gg\neq 0$ we have that
$$\bar u(r)\sim \gg\ln\frac 1r\quad\text{as }r\to 0.
$$
If we use again $\overline{e^u}\geq e^{\bar u}$ we see that  the relation $e^u\in L^1(B_{\gr_0})$ implies $\gg\leq 2$.
\qeda\medskip\\

\nind\Remark If $N=2$ the singular solution of the eikonal equation 
  \begin{equation}\label {EIK}
m|\nabla u|^q= e^{u}\quad\text{in }\BBR^2\setminus\{0\}
 \end{equation}
is expressed by
  \begin{equation}\label {EIK1}
u_e(x)=q\ln\left(\frac{2m^\frac 1q}{|x|}\right).
 \end{equation}
 It is an harmonic function. If $1<q<2$, $\nabla u_e\in L^1_{loc}(\BBR^2)$, hence $u_e$ satisfies
also
   \begin{equation}\label {E3-XY}
-\Delta u+m|\nabla u|^q= e^{u}+2\gp q\gd_0\quad\text{in }\CD'(\BBR^2).
 \end{equation}



In the next result we consider signed solutions and we give a key two-side estimate which will prove to be fundamental in the study of  their behaviour near $0$.

\bprop {DM1} Let $N\geq 3$, $1<q<2$ and $u\in C^2(\overline B_1\setminus\{0\})$ be a solution of $(\ref{Na-1})$ such that $|x|^2e^u\in L^\infty(B_1)$. Then $|\nabla u|^q$ is integrable and there exists 
$\gg\leq 0$ such that 
\bel{S5}
-\Gd u+ m|\nabla u|^q-e^u=c_N\gg\gd_0\quad\text{in }\CD'(B_1).
\ee
Furthermore,  there holds
\bel{S6}\BA{lll}\dsps
\gg (|x|^{2-N}-1)+\min_{\prt B_1}u -m\BBG_{B_1}[|\nabla u|^q](x)\leq  u(x)\\[2mm]\phantom{---------}
\dsps\leq\min\left\{ 2\ln\frac1{ |x|}+K_1,\gg (|x|^{2-N}-1)+\frac{e^{K_1}}{N-2}\ln\frac1{ |x|}\right\}\quad\text{in } B_{1},
\EA\ee
where $\BBG_{B_1}[ \; \!.  \;\!]$ is the Green operator in $B_1$. Finally $\gg=0$ if $q\geq \frac N{N-1}$.
\es
\Proof Since $u(x)\leq \ln\frac 1{|x|^2}+K_1$ for some constant $K_1$, we define the nonnegative function $V$ by
$$V(x)=-u(x)+\ln\frac 1{|x|^2}+K_1.
$$
It satisfies
\bel{S7}
-\Gd V-\frac{2(N-2)}{|x|^2}+\frac{1}{|x|^2}e^{K_1-V}=m\left|\frac{2}{|x|}{\bf e}+\nabla V\right|^q
\ee
where ${\bf e}={\bf e}_x=\frac x{|x|}$. Since $N\geq 3$ and $V\geq 0$, both $\frac{2(N-2)}{|x|^2}$ and $\frac{1}{|x|^2}e^{K_1-V}$ are integrable. By the Brezis-Lions Lemma, $\left|\frac{2}{|x|}{\bf e}+\nabla V\right|^q$ is integrable in $B_1$ and there exists $\tilde \gg\geq 0$ such that 
\bel{S8}
-\Gd V-\frac{2(N-2)}{|x|^2}+\frac{1}{|x|^2}e^{K_1-V}=m\left|\frac{2}{|x|}{\bf e}+\nabla V\right|^q+c_N\tilde \gg\gd_0\quad\text{in }\CD'(B_1).
\ee
Replacing $V$ by $-u(x)+\ln\frac 1{|x|^2}+K_1$ it infers that $|\nabla u|^q$ is integrable and $(\ref{S5})$ holds with $\gg=-\tilde\gg$.  From $(\ref{S5})$, 
$$-\Gd u\leq e^u+ c_N\gg\gd_0.
$$
Since by assumption $u\leq K_1$ on $\prt B_1$, we have that $u\leq \gg (|x|^{2-N}-1)+u_2+K_1$ where 
$$u_2(x)=\frac{e^{K_1}}{N-2}\ln \frac{1}{|x|}.
$$
For the estimate from below we have 
$$-\Gd u\geq-m|\nabla u|^q+c_N\gg\gd_0.
$$ 
This yields the left-hand side of $(\ref{S6})$. Because $e^u$ and $|\nabla u|^q$ are integrable in $B_1$, there holds when $\gg<0$, as in the proof of \rprop{Dirac2},
\bel{S8*} \bar u_r(r) =(2-N)\gg r^{1-N}(1+o(1))\quad\text{as }r\to 0,
\ee
and $\bar u_r>0$ near $0$. Hence 
$$\frac{|\gg|^q}{2^q} r^{(1-N)q}\leq\bar  u^q_r(r)\leq \overline {(u_r^2+r^{-2}|\nabla'u|^2)^{\frac q2}}.
$$
If $q\geq \frac N{N-1}$, we obtain a contradiction with $|\nabla u|^q\in L^1(B_1)$. Note that if $1<q<\frac N{N-1}$  we have from $(\ref{S8*})$
\bel{S9}
\bar u(x)=\gg|x|^{2-N}+o(|x|^{2-N}) \quad\text{as } x\to 0.
\ee
\qeda
\medskip

\nind\Remark Using standard results on elliptic equations with measure data, we also obtain that $u\in M^{\frac N{N-2}}(B_1)$ and 
$\nabla u\in M^{\frac N{N-1}}(B_1)$, where $M^p(B_1)$ denotes the Marcinkiewicz space coinciding with the Lorentz space $L^{p,\infty}(B_1)$. As consequence, if $1<q<\frac N{N-1}$ we obtain that $\nabla u\in L^{q(1+\ge)}(B_1)$ for some $\ge>0$.\medskip
\bcor{Cor1*} Under the assumptions of \rcor{Cor1}, $\gg=0$.\es
 \Proof If we apply the Green representation to formula  $(\ref{Na-2-13})$ in \rcor {Cor1} we have
\bel{S9-0}
u(x)=\gg|x|^{2-N}+\BBG_{B_{\gr_0}}[e^{m}+G-m|\nabla u|^q]+ u_{\prt B_{\gr_0}},
\ee
which implies
$$
\bar u(r)=\gg (r^{2-N}-\gr_0^{2-N})+\overline{\BBG_{B_{\gr_0}}[e^{m}+G-m|\nabla u|^q]}(r)+ \bar u(\gr_0),.
$$
Since $e^{m}+G-m|\nabla u|^q$ is integrable in $B_{\gr_0}$ there holds
$$
\lim_{r\to 0}r^{N-2}\overline{\BBG_{B_{\gr_0}}[e^{m}+G-m|\nabla u|^q]}(r)=0.
$$
which implies $\dsps\lim_{r\to 0}r^{N-2}\bar u(r)=\gg$. Because $\overline {e^u}(r)\geq e^{\bar u(r)}$ and $e^u$ is integrable we conclude that $\gg=0$.\qeda\medskip

An unexpected consequence of \rprop{DM1} is the following weak isotropy result (compare with \rprop{int-3} in next section).

\bcor {DM1cor} Let $N\geq 3$, $q>1$ and $u\in C^2(B_1\setminus\{0\})$ be a solution of $(\ref{Na-1})$ such that $|x|^2e^u\in L^\infty_{loc}(B_1)$. Then  
\bel{S9-1}
\max_{0<|x|=|y|=r}|u(x)-u(y)|=o(r^{2-N})\quad\text{ as } r\to 0.
\ee
If we assume moreover that $|\nabla u|^q\in L^{1+a}(B_1)$ for some $a>0$, then 
\bel{S9-1*}
\max_{0<|x|=|y|=r}|u(x)-u(y)|=o(r^{2-\frac N{1+a}})\quad\text{ as } r\to 0.
\ee
\es
\Proof Averaging $(\ref{S6})$ over $S^{N-1}$ we have
$$
\gg r^{2-N} -m\overline{\BBG_{B_1}[|\nabla u|^q]}(r)\leq  \bar u(r)\leq\min\left\{ 2\ln\frac1{ r}+K_1,\gg r^{2-N}+\frac{e^{K_1}}{N-2}\ln\frac1{ r}\right\},
$$
hence, since $\gg\leq 0$,
$$\BA{lll}\dsps
u(r,\gs)-\bar u(r)\leq \min\left\{ 2\ln\frac1{ r}+K_1,\gg r^{2-N}+\frac{e^{K_1}}{N-2}\ln\frac1{ r}\right\}-\gg r^{2-N} +m\overline{\BBG_{B_1}[|\nabla u|^q]}(r)\\[2mm]
\phantom{u(r,\gs)-\bar u(r)}\dsps\leq \min\left\{ 2\ln\frac1{ r}-\gg r^{2-N}+K_1,\frac{e^{K_1}}{N-2}\ln\frac1{ r}\right\}+m\overline{\BBG_{B_1}[|\nabla u|^q]}(r)\\[2mm]
\phantom{u(r,\gs)-\bar u(r)}\dsps\leq \min\left\{2,\frac{e^{K_1}}{N-2}\right\}\ln\frac 1r+m\overline{\BBG_{B_1}[|\nabla u|^q]}(r).
\EA
$$
Similarly
$$\BA{lll}\dsps
u(r,\gs)-\bar u(r)\geq \gg r^{2-N} -m\overline{\BBG_{B_1}[|\nabla u|^q]}(r)-\min\left\{ 2\ln\frac1{ r}+K_1,\gg r^{2-N}+\frac{e^{K_1}}{N-2}\ln\frac1{ r}\right\}-\\[2mm]
\phantom{u(r,\gs)-\bar u(r)}\dsps\geq \max\left\{\gg r^{2-N}-2\ln\frac1{ r}-K_1,-\frac{e^{K_1}}{N-2}\ln\frac1{ r}\right\}-m\overline{\BBG_{B_1}[|\nabla u|^q]}(r)\\[2mm]
\phantom{u(r,\gs)-\bar u(r)}\dsps\geq -\min\left\{2,\frac{e^{K_1}}{N-2}\right\}\ln\frac 1r-m\overline{\BBG_{B_1}[|\nabla u|^q]}(r).
\EA
$$
Therefore
\bel{S9-2}
|u(r,\gs)-\bar u(r)|\leq C\ln\frac 1r+o(r^{2-N})=o(r^{2-N})\quad\text{ as }\;r\to 0.
\ee
 If  $|\nabla u|^q\in L^{1+a}(B_1)$, then 
$$\BA{lll}\dsps 0\leq \overline{\BBG_{B_1}[|\nabla u|^q]}(r)=\int_0^rs^{1-N}\int_0^s\overline{|\nabla u|^q}t^{N-1}dtds\\[4mm]
\phantom{\dsps 0\leq \overline{\BBG_{B_1}[|\nabla u|^q]}(r)}\dsps\leq \int_0^rs^{1-N}\left(\int_0^s\overline{|\nabla u|^{q(1+a)}}t^{N-1}dt\right)^{\frac{1}{1+a}}
\left(\int_0^st^{N-1}dt\right)^{\frac{a}{1+a}}.
\EA$$
This implies $(\ref{S9-1*})$.

\qeda
\medskip

\nind\Remark When $|x|^2e^u$ remains bounded, then for any $q>1$  we have that $|\nabla u(x)|^q\leq c|x|^{-\frac{q}{q-1}}$ by \rcor{Th1}; hence  in the supercritical case $q>\frac N{N-1}$ we have that
$|\nabla u|^q\in L^{1+a}(B_1)$ for any $0<a<\frac{N(q-1)-q}{q}$. 
\medskip

Thanks to the unconditional estimate of \rth{Th3} we have an extension of \rprop{DM1} to the case $N>q>2$. 

\bprop {DM2} Let $N>q>2$ and $u\in C^2(\overline B_1\setminus\{0\})$ be a solution of $(\ref{Na-1})$. Then $e^u$ and $|\nabla u|^q$ are integrable and $u$ satisfies
\bel{S10}
-\Gd u+ m|\nabla u|^q-e^u=0\quad\text{in }\CD'(B_1).
\ee
Furthermore,  there holds for some $K>0$,
\bel{S11*}
-m\BBG_{B_1}[|\nabla u|^q](x)+\min_{\prt B_1}u\leq  u(x)\leq q\ln\frac1{ |x|}+K\quad\text{in } B_{1}.
\ee
\es
\Proof We mimic the proof of \rprop{DM1}. By \rth{Th3} there exists a real number $K_1$ such that $u(x)\leq q\ln\frac1{|x|}+K_1$, therefore the function
 $V(x)=-u(x)+q\ln\frac{1}{ |x|}+K_1$ is nonnegative and satisfies
\bel{S12}
-\Gd V-\frac{q(N-2)}{|x|^2}+\frac{1}{|x|^q}e^{K_1-V}=m\left|\frac{q}{|x|}{\bf e}+\nabla V\right|^q.
\ee
Since $q<N$, $\frac{1}{|x|^q}e^{K_1-V}$ is integrable and clearly $\frac{q(N-2)}{|x|^2}$ is so. Then by the Brezis-Lions Lemma $\left|\frac{q}{|x|}{\bf e}+\nabla V\right|^q$ is integrable and there exists $\tilde \gg\geq 0$ such that 
\bel{S13}
-\Gd V-\frac{q(N-2)}{|x|^2}+\frac{1}{|x|^q}e^{K_1-V}=m\left|\frac{q}{|x|}{\bf e}+\nabla V\right|^q+c_N\tilde \gg\gd_0\quad\text{in }\CD'(B_1).
\ee
Since $q<N$, $|x|^{-q}$ is integrable, thus $[\nabla V|^q$ is integrable too and $[\nabla u|^q$ shares this property. If we had $\gg=-\tilde \gg<0$ and since $(\ref{S8*})$ holds we would  reach a contradiction since is is assumed that $q\geq \frac N{N-1}$. As a consequence we have that $\gg=0$, hence  $\dsps\lim_{r\to 0}r^{N-2}\bar u(r)=0$ and therefore estimate $(\ref{S6})$ turns into 
$(\ref{S11*})$.\qeda\medskip

\mysection{Isolated singularities in the case $1<q<2$}

This section is devoted to proving  \rth{Th5}. Since it is long and involves very different techniques, we divide our proof into several steps.
\subsection{Some tools}

In \cite{BV-V1} the main tool for studying positive solutions was based upon Simon's fundamental work on analytic functional. Because of the presence of the term 
$|\nabla u|^q$  we use the following extension due to Huang and Tak\'a$\rm\check{c}$ \cite{HuTak} of Simon's results. 
In their framework we denote by  $(M,g)$ a compact Riemannian manifold with tangent bundle $TM$,  with norm $|.|_g$, and by $\Gd_g$ the Laplace-Beltrami operator on $M$.
\bth{th**} Let  $f:\BBR\mapsto\BBR$ be a real analytic function and $g:\BBR_+\ti\BBR\ti TM$ be a $C^1$ function such that 
\bel{SN-1}
|g(t,r,\xi)|\leq \eta(t)a(r)b(|\xi|_g)
\ee
where $a$ and $b$ are locally bounded functions and $\eta$ satisfies for some $\gg>0$,
\bel{SN2}
\int_t^\infty\eta^2(s)ds\leq ct^{-1-\gg}\quad\text{for all }t\geq 1.
\ee
If $u$ is bounded on $[0,\infty)\ti M$ and satisfies 
\bel{SN3}
u_{tt}+\ge u_t+\Gd_gu+f(u)+g(t,u,\nabla_gu)=0\quad\text{in }[0,\infty)\ti M,
\ee
where $\ge\neq 0$, there there exists $\gw\in C^2(M)$ such that 
\bel{SN4}
\Gd_g\gw+f(\gw)=0\quad\text{in }M,
\ee
such that $u(t,.)$ converges to $\gw$ in $C^2(M)$ when $t\to\infty$
\es

The second tool is the isotropy theorem, already  used in a weaker form in \cite{BV-V1}, which actually relies on a representation formula adapted from \cite[Appendix]{Ve1} (see also \cite[Theorem 1.1]{Bo-Ve}).
We denote by $(r,\gs)\in \BBR_+\ti S^{N-1}$ the spherical coordinates in $\BBR^N$. If $u$ is a solution of $(\ref {Na-1})$  in $B_1\setminus\{0\}$ it satisfies
\bel{SN5}u_{rr}+\frac {N-1}{r}u_r+\frac1{r^2}\Gd_gu+e^{u}-m\left(u_r^2+r^{-2}|\nabla_gu|^2\right)^{\frac q2}=0\quad\text{in }(0,1)\ti S^{N-1}.
\ee
We transform it into an equation in the cylinder $(0,\infty)\ti S^{N-1}$ by setting
\bel{SN6}r=e^{-t}\Longleftrightarrow t=\ln\frac 1{r}\quad\text{and }\,v(t,\gs)=u(r,\gs)-2\ln\frac 1{r}\Longrightarrow v(t,\gs)=u(r,\gs)-2t.\ee
Then there holds in $(0,\infty)\ti S^{N-1}$
\bel{SN7}v_{tt}-(N-2)v_t+\Gd_gv-2(N-2)+e^{v}-me^{(q-2)t}\left(\left(v_t+2\right)^2+|\nabla_gv|^2\right)^{\frac q2}=0.
\ee
If $u$ is a solution of $(\ref{Na-1})$ in $B_1\setminus\{0\}$ such that $|x|^2e^u\leq\tilde K_1$, then the function $v$ satisfies 
\bel{SN8} v(t,\gs)\leq  K_1.
\ee
In the next proposition we give a regularity estimate on $v$ provided it satisfies some estimate from below.

\blemma{int-2} Let $1<q<2$ and $u$ be a solution $(\ref{Na-1})$ in $B_1\setminus\{0\}$ such that the function $v$ defined by $(\ref{SN6})$ satisfies 
\bel{SN9} 
-\eta(t)\leq v(t,\gs)\leq K_1\quad\text{for }t\geq 1,
\ee
for some positive function $\eta\in C^2([0,\infty))$ verifying 
\bel{SN10} 
\frac{1}{\eta(t)}+\left|\frac{\eta_{t}(t)}{\eta(t)}\right|+\left|\frac{\eta_{tt}(t)}{\eta(t)}\right|\leq C\quad\text{and }\;\eta(t)\leq Ce^{\gb t}\quad\text{for }t\geq 1,
\ee
where $\gb:=\frac{2-q}{q-1}$. Then 
\bel{SN11} 
\left|v_t(t,\gs)\right|+\left|\nabla_gv(t,\gs)\right|\leq C\eta(t)\quad\text{for }(t,\gs)\in [2,\infty)\ti S^{N-1}.
\ee
Equivalently, if $|x|^2e^{u}\in L^\infty_{loc}(B_1)$ and $u(x)\geq -\eta(\ln|x|)$, then we have
\bel{SN12} 
\left|\nabla u(x)\right|\leq C\frac{\eta(-\ln|x|)}{|x|}\quad\text{for }x\in B_{e^{-2}}\setminus\{0\}.
\ee
In the particular case where $\eta (t)=e^{\gn t}$ with $\gn\leq \gb$ the previous result asserts that 
\bel{SN14} 
-K_1|x|^{-\gn}\leq u(x)\leq 2\ln \frac1{|x|}+K_2\Longrightarrow |\nabla u(x)|\leq K_3|x|^{-\gn-1}\quad\text{for all  }x\in B_{\frac12}\setminus\{0\}.
\ee
\es
\Proof Set $v(t,\gs)=\eta(t)\psi(t,\gs)$. Then
$$\BA{lll}\dsps \psi_{tt}-\left(N-2-2\frac{\eta_{t}}{\eta}\right)\psi_t+\left(\frac{\eta_{tt}}{\eta}+(2-N)\frac{\eta_{t}}{\eta}\right)\psi+\Gd_g\psi\\[4mm]\phantom{--------}\dsps
-\frac{2(N-2)}{\eta}+\frac{e^{\eta\psi}}{\eta}-me^{(q-2)t}\eta^{q-1}\left(\left(\psi_t+\frac{\eta_{t}}{\eta}\psi-\frac{2}{\eta}\right)^2+|\nabla'\psi|^2\right)^{\frac q2}=0.
\EA$$
Since $\eta^{-1}$ is bounded and $(\ref{SN10} )$ holds, we can apply the standard elliptic equations results and derive that $\psi_t$ and $|\nabla'\psi|$ are uniformly bounded. Hence $(\ref{SN11})$ and $(\ref{SN12} )$ hold.\qeda

\medskip

\bprop{int-3} Let $N\geq 3$, $1<q<2$ and $u$ be a solution $(\ref{Na-1})$ in $B_1\setminus\{0\}$ such that the function $v$ defined by $(\ref{SN6})$ satisfies 
\bel{SN15} 
-Ce^{\gn t}\leq v(t,\gs)\leq K_1\quad\text{for }t\geq 1,
\ee
for some $0\leq \gn\leq \gb$ satisfying also $\gn<\min\left\{N-1,\frac{N+1}{q}-1\right\}$, then 
\bel{SN16} \dsps
|u(x)-u(y)|\leq K_5r^{2-q-q\gn}+K_6\quad\text{for all }0<|x|=|y|=r\leq\frac 12, 
\ee
for some $K_5,K_6>0$.
\es
\Proof {\it Step 1: Equation in cylinder}. Set $v(t,\gs)=e^{\gn t}\psi(t,\gs)$. Then $\psi$ is bounded by assumption and it satisfies
$$\BA{lll}\dsps \psi_{tt}-\left(N-2-2\gn\right)\psi_t+\gn(\gn+2-N)\psi+\Gd_g\psi-2(N-2)e^{-\gn t}+e^{-\gn t}e^{e^{\gn t}\psi}\\[4mm]\phantom{-------------}\dsps
 -me^{((q-2)+(q-1)\gn)t}\left(\left(\psi_t+\gn\psi-2e^{-\gn t}\right)^2+|\nabla_g\psi|^2\right)^{\frac q2}=0.
\EA$$
If $\psi^*(t,.)=\psi(t,.)-\overline {\psi}(t)$, then 
\bel{SN17}\BA{lll}\dsps \psi^*_{tt}+a\psi^*_t-\ell\psi^*+\Gd_g\psi^*+F^*(t)=0,
\EA\ee
where 
\bel{SN18^}\ell=\gn(N-2-\gn)\text{ and }\;\;a=2\gn+2-N,\ee
and
\bel{SN18}\BA{lll}F^*(t)=e^{-\gn t}\left(e^{e^{\gn t}\psi}-\overline {e^{e^{\gn t}\psi}}\right)\\[2mm]\phantom{-}-me^{((q-2)+(q-1)\gn)t}\left(\left(\left(\psi_t+\gn\psi-2e^{-\gn t}\right)^2+|\nabla_g\psi|^2\right)^{\frac q2}-\overline {\left(\left(\psi_t+\gn\psi-2e^{-\gn t}\right)^2+|\nabla_g\psi|^2\right)^{\frac q2}}\right)\\
\phantom{F^*(t}=F^*_1(t)+F^*_2(t),\EA\ee
in which expression we have set
$$F^*_1(t)=e^{-\gn t}\left(e^{e^{\gn t}\psi}-\overline {e^{e^{\gn t}\psi}}\right)\quad\text{and }\;F^*_2(t)=F^*(t)-F^*_1(t).$$
By \rlemma{int-2} and the fact that $\gn\leq \gb$ the functions $F_j^*$ ($j=1,2$) satisfy
\bel{SN18*}\BA{lll}
(i)\qquad\qquad&\norm {F^*_1(t)}_{L^\infty}\leq e^{-\gn t}\qquad\qquad\qquad\qquad\qquad\qquad\qquad\qquad\qquad\qquad\qquad\quad\quad\\[1mm]
(ii)\qquad\qquad&\norm {F^*_2(t)}_{L^\infty}\leq c_0e^{(q-2+(q-1)\gn)t}.
\EA\ee

\nind {\it Step 2: The representation formula}. 
We consider the fractional operator  $\BBL$,
$$\BBL=-\left(-\Gd_g+\left(\tfrac {a^2}4+\ell\right)I\right)^{\frac 12}=-\left(-\Gd_g+\tfrac {(N-2)^2}4I\right)^{\frac 12}.
$$
It is maximal monotone in $L^2(S^{N-1})$ with domain $W^{1,2}(S^{N-1})$. For any function $\gf\in W^{2,2}(S^{N-1})$ with zero average there holds
$$\int_{S^{N-1}}\phi\left(-\Gd_g\phi+\tfrac {(N-2)^2}{4}\phi\right)dS\geq \left(N-1+\tfrac {(N-2)^2}4\right)\int_{S^{N-1}}\phi^2dS=\frac {N^2}4\int_{S^{N-1}}\phi^2dS.
$$
Then the restriction of the semigroup $S(t)=e^{t\BBL}$ to the subspace $H$ of $L^2(S^{N-1})$ of functions with zero average satisfies
\bel{SN19}
\norm{ S(t)[\gf]}_{L^2(S^{N-1})}\leq e^{-\frac N2t}\norm\phi_{L^2(S^{N-1})}\quad\text{for all }\phi\in H.
\ee
 By standard Fourier analysis, the restriction of $S(t)$ to $H\cap L^\infty(S^{N-1})$ satisfies
 \bel{SN20}
\norm{ S(t)[\gf]}_{L^\infty(S^{N-1})}\leq Ce^{-\frac N2t}\norm\phi_{L^\infty(S^{N-1})}\quad\text{for all }\phi\in H\cap L^\infty(S^{N-1}),
\ee
for some $C=C(N)>0$. Next we set for $t\geq 0$, 
\bel{SN23}\BA{lll}\dsps\Psi^*(t)=e^{-\frac{at}2}S(t)[\psi^*(0)]+\int_0^te^{-\frac{as}2}S(s)\left[\int_0^\infty e^{\frac{a\gt}2}S(\gt)
\left[F^*(t+\gt-s)\right]d\gt\right]ds.
\EA\ee
We first check that the integral term on the right-hand side is normally convergent.
From $(\ref{SN20})$ there holds
$$\norm{e^{\frac{at}{2}}S(t)}_{\CL(H\cap L^\infty,H\cap L^\infty)}\leq c_3e^{(\gn+1-N)t},
$$
and
$$\norm{e^{-\frac{at}{2}}S(t)}_{\CL(H\cap L^\infty,H\cap L^\infty)}\leq c_4e^{-(\gn+1)t}.
$$
Formally, 
$$\int_0^\infty e^{\frac{a\gt}2}S(\gt)
\left[F^*(t+\gt-s)\right]d\gt=\int_0^\infty e^{\frac{a\gt}2}S(\gt)
\left[F_1^*(t+\gt-s)\right]d\gt+\int_0^\infty e^{\frac{a\gt}2}S(\gt)
\left[F_2^*(t+\gt-s)\right]d\gt.
$$
Since $\gn<N-1$ and using $(\ref{SN18*})$-(i) the integral over $(0,\infty)$ of the term involving $F^*_1$ is normally convergent and we have
\bel{SN24}\norm{\int_0^\infty e^{\frac{a\gt}2}S(\gt)
\left[F_1^*(t+\gt-s)\right]d\gt}_{L^\infty}\leq c_3e^{-\gn(t-s)}.
\ee
We have also from  $(\ref{SN18*})$-(ii) 
$$\norm{e^{\frac{a\gt}2}S(\gt)
\left[F_2^*(t+\gt-s)\right]}_{L^\infty}\leq c_0c_3e^{(q-2+\gn(q-1))(t-s)}e^{((\gn+1)q-N-1)\gt}.
$$
Hence, the convergence of the integral involving $F_2^*(t+\gt-s)$ is guaranteed if $\gn<\frac{N+1-q}{q}$, 
and we obtain
\bel{SN24*}\norm{\int_0^\infty e^{\frac{a\gt}2}S(\gt)
\left[F_2^*(t+\gt-s)\right]d\gt}_{L^\infty}\leq c_4e^{(q-2+(q-1)\gn)(t-s)}.
\ee
Furthermore we obtain the following estimates
\bel{SN24**}\norm{\int_0^te^{-\frac {as}2}S(s)\left[\int_0^\infty e^{\frac{a\gt}2}S(\gt)
\left[F_1^*(t+\gt-s)\right]d\gt\right]
}_{L^\infty}\leq c_5e^{-\gn t},
\ee
and
\bel{SN24***}\norm{\int_0^te^{-\frac {as}2}S(s)\left[\int_0^\infty e^{\frac{a\gt}2}S(\gt)
\left[F_2^*(t+\gt-s)\right]d\gt\right]
}_{L^\infty}\leq c_6e^{(q-2+(q-1)\gn)t}.
\ee
As a consequence
\bel{SN25}\BA{lll}\dsps
\norm{\Psi^*(t)}_{L^\infty}\leq c_7e^{-(\gn+1)t}+c_5e^{-\gn t}+c_6e^{(q-2+(q-1)\gn)t}\leq c_8e^{-\gn t}+c_6e^{(q-2+(q-1)\gn)t}.
\EA\ee

\nind{\it Step 3: End of the proof}. We have now to prove that $\psi^*=\Psi^*$. Clearly $\Psi^*(0)=\psi^*(0)$. Then we put
\bel{SN23+}\BA{lll}\dsps\Phi(t)=\int_0^te^{-\frac{as}2}S(s)\left[\int_0^\infty e^{\frac{a\gt}2}S(\gt)
\left[F^*(t+\gt-s)\right]d\gt\right]ds\\[4mm]
\phantom{\dsps\Phi(t)}\dsps=e^{-\frac{at}2}\int_0^tS(t-s)\left[\int_s^\infty e^{\frac{a\gt}2}S(\gt-s)\left[F^*(\gt)\right]d\gt\right]ds.
\EA\ee
Since $\BBL$ is a closed operator and $\Phi$ has all the needed regularity requirements, we have
$$\Phi_t(t)=(\BBL-\frac a2I)\Phi(t)+e^{-\frac{at}{2}}\int_t^\infty e^{\frac{a\gt}2}S(\gt-t)\left[F^*(\gt)\right]d\gt,
$$
$$\BA{lll}\dsps\Phi_{tt}(t)=(\BBL-\frac a2I)\Phi_t(t)-e^{-\frac{at}{2}}(\BBL+\frac a2I)\int_t^\infty e^{\frac{a\gt}2}S(\gt-t)\left[F^*(\gt)\right]d\gt-F^*(t)\\[4mm]
\phantom{\dsps\Phi_{tt}(t)}\dsps=(\BBL-\frac a2I)^2\Phi(t)-ae^{-\frac{at}{2}}\int_t^\infty e^{\frac{a\gt}2}S(\gt-t)\left[F^*(\gt)\right]d\gt-F^*(t).
\EA$$
Then
$$\BA{lll}\dsps\CE(\Phi)(t):=\Phi_{tt}(t)+a\Phi_t(t)-\ell\Phi(t)+\Gd_g\Phi(t)+F^*(t)\\[4mm]
\phantom{\dsps\CE(\Phi)(t):}\dsps=(\BBL-\frac a2I)^2\Phi(t)-ae^{-\frac{at}{2}}\int_t^\infty e^{\frac{a\gt}2}S(\gt-t)\left[F^*(\gt)\right]d\gt\\[4mm]
\phantom{\dsps\CE(\Phi)(t):=}\dsps+a(\BBL-\frac a2I)\Phi(t)+ae^{-\frac{at}{2}}\int_t^\infty e^{\frac{a\gt}2}S(\gt-t)\left[F^*(\gt)\right]d\gt+(\Gd_g-\ell I) \Phi(t)\\[4mm]
\phantom{\dsps\CE(\Phi)(t):}\dsps= \left(\BBL^2-\left(-\Gd_g+\left(\tfrac{a^2}{4}+\ell\right)I\right)\right) \Phi(t)=0,
\EA$$
from the definition of $\BBL$. \\
Performing the same computation as above it is easy to verify that the function $(t,\gs)\mapsto e^{-\frac{at}2}S(t)[\psi^*(0)](\gs)$ is a solution of the mere equation
 \bel{SN21} y_{tt}+ay_t-\ell y+\Gd_gy=0\quad\text{on }\BBR^+\ti S^{N-1},
\ee
with zero spherical average as $\psi^*$, and we have proved that  $\Phi$ satisfies it. Since $\Psi^*(t)=\Phi(t)+e^{-\frac{at}2}S(t)[\psi^*(0)]$, the function $(t,.)\mapsto\Psi^*(t,.)-\psi^*(t,.)$ satisfies $\ref{SN21})$ on $\BBR^+\ti S^{N-1}$ and vanishes at $t=0$. From the assumption and relation $(\ref{SN25})$ we have that it is  is bounded. Furthermore, since it has zero average, the function $t\mapsto Z(t):=\norm{\Psi^*(t,.)-\psi^*(t,.)}_{L^2(S^{N-1})}$ satisfies
 \bel{SN21*} 
 Z_{tt}+aZ_t-(\ell+N-1)Z\geq 0
 \quad\text{in }\CD'(\BBR^+).
\ee
The characteristic roots of the equation
 \bel{SN21*$} 
z_{tt}+az_t-(\ell+N-1)z= 0
\ee
are $\gl_+=\frac{-a+\sqrt{a^2+4(\ell+N-1)}}{2}=N-1-\gn$ and $\gl_-=-1-\gn$, using relations $(\ref{SN18^})$. Notice that $\gl_+>0$ and 
$\ell+N-1=(\gn+1)(N-1-\gn)>0$ by the assumption on $\gn$. Since $Z(0)=0$ and $Z(t)=o(e^{\gl_+t})$ as $t\to\infty$, we have that $Z\equiv 0$, hence 
$\Psi^*=\psi^*$ and $\psi^*$ satisfies the relation $(\ref{SN23})$ as $\Psi^*$ does. This implies that estimate $(\ref{SN25})$ is verified by $\psi^*$. 
and therefore
\bel{SN26}
\norm{v(t,.)-\overline {v}(t)}_{L^\infty}\leq c_6e^{(\gn q+q-2)t}+c_8,
\ee
from which $(\ref{SN16})$ follows.\qeda\medskip

\nind\Remark Formula $(\ref{SN23})$ shows that $\Psi^*$ (and thus $\psi^*$) is solution of a system of differential equations in the space $\BBR_+\ti H$
\bel{Sys1}
\left\{\BA{lll}\frac{d\Psi^*}{dt}+\BBL\Psi^*+\frac a2\Psi^*=\Xi\\[2mm]
\phantom{--}\frac{d\Xi}{dt}-\BBL\Xi+\frac a2\Xi=F^*\\[2mm]
\Psi^*(0)=\psi^*(0)\,\text{ and }\dsps\lim_{t\to\infty}\Xi(t)=0.\EA
\right.
\ee
The decomposition of second order elliptic equations in cylinders is originally due to S.G. Krein \cite{Kr}. It was rediscovered and used in \cite{Ve1} for the study of isolated singularities of solutions of the Emden-Fowler equation 
\bel{Sys2}
-\Gd u+|u|^{q-1}u=0.
\ee
Its was also used under a simpler form in order to obtain isotropy estimates in \cite{BV-V1}. In \cite{Bo-Ve} this method allows to obtain a representation theorem valid for any positive solutions of 
\bel{Sys3}
v_{tt}+av_t-\ell v+\BBA v+v^q=0
\ee
on a Riemannian manifold $(M,g)$ where $\BBA$ is a uniform elliptic operator on $M$. Several uniform estimates follow from this representation.\medskip

\nind\Remark If $|x|^2e^{u}\in L^\infty(B_1)$, then $|x|^{\gb}u(x)$ is uniformly bounded by \rcor{Th1}. Thus $\gn=\gb$. The assumption $\gn<\frac{N+1}{q}-1$ is satisfied 
if $q>\frac {N+1}N$, but the estimate $(\ref{SN16})$ brings nothing new because $2-q-q\gn=-\gb$. \medskip

\nind{\it Proof of \rcor{Harn1} i.e. Harnack inequality}. Note that the assumption contains the two cases $u(x)-a>b$ and $u(x)-a<-b$. We set $u_a=u-a$. Then $u_a$ satisfies
$$-\Gd u_a+m|\nabla u_a|^q-e^ae^{u_a}=0.
$$
We write $|\nabla u_a|^q=|\nabla u_a|^{q-1}|\nabla u_a|$ and $e^ae^{u_a}=\frac{e^ae^{u_a}}{u_a}u_a$. By $(\ref{Na-17})$ and $(\ref{Na-18})$, $|x|^2\frac{e^ae^{u_a}}{u_a}$ and  $|x||\nabla u_a|^{q-1}|$ are bounded near $x=0$. Therefore (see e.g. \cite{GiTr}) $u_a$ satisfies Harnack inequality in the sense that  $(\ref{Na-20})$ holds.\qeda

\subsection{Positive singularities}

Our main result relies directly on   \rth{th**} and \rlemma{int-2}.

\bprop{conv} Let  $N\geq 3$, $1<q<2$ and  $u$ be a solution of $(\ref{Na-1})$ in $B_1\setminus\{0\}$ such that  
\bel{SN30*}
0<\liminf_{x\to 0}|x|^2e^{u(x)}\leq \limsup_{x\to 0}|x|^2e^{u(x)}<\infty.
\ee
Then there exists a solution $\gw$ of $(\ref{Na-14})$ such that 
\bel{SN28}
u(r,\gs)-2\ln\frac 1r\to \gw(\gs)\quad\text{as } r\to 0,
\ee
uniformly on $S^{N-1}$.
\es
\Proof The assumption $(\ref{SN30*})$ implies that $u(x)-2\ln\frac 1{|x|}$ is bounded in $B_1$. Then the function $v$ defined by $(\ref{SN6})$ is bounded. By \rlemma{int-2} the functions $v_t$,  
$\nabla_gv$, $v_{tt}$, $\nabla v_t$ and $D^2v$  remain also uniformly bounded. From $(\ref{SN7})$ there holds
$$\BA{lll}\dsps\frac{d}{dt}\int_{S^{N-1}}\left(\frac{1}{2}\left(v_t^2-|\nabla_gv|^2\right)-2(N-2)v+e^v\right) dS\\[4mm]
\phantom{----}\dsps=(N-2)\int_{S^{N-1}}v_t^2dS+me^{(q-2)t}\int_{S^{N-1}}\left((v_t+2)^2+|\nabla_gv|^2\right)^{\frac q2}v_tdS.
\EA$$
The energy function
$$\CE[v](t):=\int_{S^{N-1}}\left(\frac{1}{2}\left(v_t^2-|\nabla_gv|^2\right)-2(N-2)v+e^v\right) dS
$$
remains uniformly bounded. Since $q<2$ we obtain that
$$\int_0^\infty\int_{S^{N-1}}v_t^2dSdt<\infty.
$$
Since $v_t$ is uniformly continuous, it follows that $v_t(t,.)\to 0$ in $L^2(S^{N-1})$ when $t\to\infty$. Differentiating the equation and multiplying by  $v_{tt}$ yields
$$\int_0^\infty\int_{S^{N-1}}v_{tt}^2dSdt<\infty.
$$
and again $\dsps\lim_{t\to\infty}\norm{v_{tt}(t,.)}_{L^2(S^{N-1})}= 0$. As a consequence, the limit set $\Gw[v]$ at infinity of the trajectory $\dsps\CT[v]:=\cup_{t>0}\{v(t,.)\}$ in $C^2(S^{N-1})$ which is defined by 
$$\dsps \Gw[v]=\bigcap_{t>0}cl_{C^2}\left(\bigcup_{\gt\geq t}\{v(t,.)\}\right)
$$
is a non-empty compact and connected subset of the set of solutions $\gw$ of $(\ref{Na-14})$. The Huang-Takac extension of Simon's result \cite{HuTak} quoted in® \rth{th**}   asserts that $\Gw[v]=\{\gw\}$ for some $\gw$ satisfying $(\ref{Na-14})$.\qeda
\subsection{Non-positive singular solutions}

\blemma{conv-2} Let $1<q<\frac N{N-1}$, $N\geq 3$. If  $u\in C^2(\overline B_1\setminus\{0\})$ is a solution of $(\ref{Na-1})$ in $B_1\setminus\{0\}$ verifying $|x|^{2}e^{u(x)}\leq C$ and  
\bel{SN29}
\liminf_{x\to 0}|x|^{2}e^{u(x)}=0.
\ee
Then the following alternative holds:\smallskip

\nind 1-Either there exists $\gg<0$ such that 
\bel{SN30}
\lim_{x\to 0}|x|^{N-2}u(x)=\gg.
\ee
and still $u$ satisfies $(\ref{S5})$ and $(\ref{S6})$.\smallskip

\nind 2- Or $u$ can be extended by continuity as a solution of  $(\ref{Na-1})$ in $B_1$.
\es
\Proof 
By \rprop{DM1}  $|\nabla u|^q$ and $e^u$ are integrable and there exists $\gg\leq 0$ such that  $(\ref{S5})$ holds. We have already seen that $u\in M^{\frac{N}{N-2}}(B_1)$ and 
$\nabla u\in M^{\frac{N}{N-1}}(B_1)\subset L^{\frac{N-\ge}{N-1}}(B_1)$ for any $\ge>0$. This implies that $|\nabla u|^{q-1}\in L^{\frac{N-\ge}{(N-1)(q-1)}}(B_1)$, and since 
$q<\frac{N}{N-1}$, we have that $\frac 1{(q-1)(N-1)}=1+\gd$ for some $\gd>0$. We chose $\ge=\frac{2N\gd}{3+2\gd}$ and we obtain that $|\nabla u|^{q-1}\in L^{N+\frac\ge 2}(B_1)$. 
We encounter two possibilities.\smallskip

\nind {\it Case 1}: Suppose that $\gg<0$. By \rprop{DM1} and $(\ref{S9})$, $\bar u(r)=\gg r^{2-N}(1+o(1))$ as $r\to 0$ and from $(\ref{S9-1})$, 
$\lim_{x\to 0}|x|^{N-2}|u(x)-\bar u(|x|)|=0$. This implies $(\ref{SN30})$. Using now  $(\ref{SN16})$ with $\gn=N-2$, we obtain actually the following estimate
\bel{SN30+}
|\gg-|x|^{N-2}u(x)|\leq K_8|x|^{N-q(N-1)}\quad\text{in }B_1\setminus\{0\}.
\ee
\smallskip

\nind {\it Case 2}: Suppose that $\gg=0$. Then we claim that $u$ is bounded. \smallskip

\nind {\it Step 1:  We prove that $u$ is bounded from below}. We set $\CU_+:=\{x\in \bar B_1\setminus\{0\}:u(x)\geq 0\}$ and $\CU_-:=\{x\in \bar B_1\setminus\{0\}:u(x)< 0\}$. The function $u$ satisfies in $\CD'(B_1)$,
\bel{SN3X}-\Gd u+m|\nabla u|^q=e^u{\bf I{_{_{\CU_+}}}}+e^u{\bf I{_{_{\CU_-}}}}:=F,
\ee
where ${\bf I{_{_{\CU_+}}}}$ and ${\bf I{_{_{\CU_-}}}}$ denote the characteristic functions of  $\CU_+$ and  $\CU_-$. By assumption
$$0\leq  F\leq \frac K{|x|^2}{\bf I{_{_{\CU_+}}}}+{\bf I{_{_{\CU_-}}}}.$$
Hence $F\in M^{\frac{N}{2}}(B_1)\subset L^p(B_1)$ for any $1<p<\frac N2$.
Also $|\nabla u|\in M^{\frac{N}{N-1}}$, hence $|\nabla u|^q\in M^{\frac{N}{q(N-1)}}$. By regularity theory (e.g. \cite{GiTr}) we deduce that $u\in W^{2,a}(B_1)$ for any 
$a<\min\left\{\frac{N}{q(N-1)},\frac N2\right\}$. \\
(i)- If $\min\left\{\frac{N}{q(N-1)},\frac N2\right\}=\frac N2$, which is the case if $q\leq \frac 2{N-1}$, by Sobolev imbedding theorem we have that 
$\nabla u\in L^{a^*}(B_1)$ for any $a^*<N$. Since $q<2$, $a^*$ can be chosen such as $\frac {a^*}{q}>\frac N2$, thus $|\nabla u|^q\in L^b(B_1)$ with $b>\frac N2$, and then $\BBG_{B_1}[|\nabla u|^q]\in L^\infty(B_1)$. Since by $(\ref{S6})$ with $\gg=0$, we have that
\bel{SN3Y}
\min_{\prt B_1}u -m\norm{\BBG_{B_1}[|\nabla u|^q]}_{L^\infty}\leq  u(x) \leq\min\left\{ 2\ln\frac1{ |x|}+K_1,\frac{e^{K_1}}{N-2}\ln\frac1{ |x|}\right\}\quad\text{in } B_{1}.
\ee
It follows that $u$ is uniformly bounded from below in $B_1\setminus\{0\}$.\\
(ii)- If $q>\frac 2{N-1}$, then $\frac{N}{q(N-1)}<\frac N2$, then $u\in W^{2,a}(B_1)$ for any $1<a<\frac{N}{q(N-1)}$  We choose $a_0$ such that $1<a_0<\frac{N}{q(N-1)}$. Then $u\in W^{2,a_0}(B_1)$. Since $a_0<N$ we have by Sobolev imbedding, $|\nabla u|\in L^{\frac{Na_0}{N-a_0}}(B_1)$, equivalently $|\nabla u|^q\in L^{\frac{Na_0}{q(N-a_0)}}(B_1):=L^{a_1}(B_1)$. This defines 
$a_1=\frac{Na_0}{q(N-a_0)}$, and $a_1>a_0$ since $a_0>1>\frac{N(q-1)}{N}$. If $a_1\geq N$ we have the result as in (i). If we assume that $a_1< N$
we assume by induction that we have constructed $a_1<a_2<...<a_{n-1}<N$. Then we set $a_n=\frac{Na_{n-1}}{q(N-a_{n-1})}$ with the property that $|\nabla u|^q\in L^{a_{n-1}}(B_1)$ implies 
$|\nabla u|^q\in L^{a_{n}}(B_1)$. By induction $\{a_n\}$ is increasing. If it remains always bounded from above by $N$, there would exists $\ell\leq N$ such that 
$\ell=\frac{N\ell}{q(N-\ell)}$, that is $\ell=\frac {N(q-1)}{q}$, which is smaller than $1$, contradiction. Hence there exists an integer $n_0$ such as $a_{n_0}\geq N$. replacing $a_0$ we obtain the existence of some $a_{n_0}\geq N$ and we end the proof as in (i) and conclude again that $\BBG_{B_1}[|\nabla u|^q]\in L^\infty(B_1)$, thus $u$ is bounded from below.
\smallskip

\nind It follows from \rprop{int-3} with $\gn$ arbitrarily small that
\bel{SN48}|u(x)-u(y))|\leq K_6\quad \text{for } |x|=|y|\leq \frac 12.
\ee
{\it Step 2: We prove that $u$ is bounded from above}. We proceed by contradiction, assuming  that $u$ is not upper bounded. Hence there exists a sequence $\{x_n\}$ converging to $0$ such that $u(x_n)\to\infty$. Set $r_n=|x_n|$; using $(\ref{SN48})$ we deduce that 
$u(r_n,.)\to\infty$ uniformly on $S^{N-1}$ when $n\to\infty$. Put $b_n=\min\{u(z):|z|=r_n\}$. We can assume that $\{b_n\}$ is increasing. Since 
$$-\Gd u+m|\nabla u|^q\geq 0\; \text{ in }\{x: r_{n+k}<|x|<r_{n}\},
$$
we obtain that $u\geq h\geq b_n$ in $\overline B_{r_n}\setminus B_{r_{n+k}}$, by comparison with the solution of 
$$\BA{lll}-\Gd h+m|\nabla h|^q= 0 &\text{in }\{x: r_n<|x|<r_{n+k}\}\\
\phantom{-\Gd +m|\nabla h|^q}
h=b_n&\text{if }|x|=r_n\\
\phantom{-\Gd +m|\nabla h|^q}
h=b_{n+k}&\text{if }|x|=r_{n+k}.
\EA$$
Letting $k\to\infty$ yields $u\geq b_n$ in $\overline B_{r_n}\setminus\{0\}$. Therefore
\bel{SN49}
\lim_{x\to 0}u(x)=\infty.
\ee
The function  $v$ defined by $u(x)=v(t,.)+2t$ with $t=\ln \frac 1r$ satisfies for some $T>0$
\bel{SN47X}
C-2t\leq v(t,\gs)\leq K_1\quad\text{in }[T,\infty)\ti S^{N-1},
\ee
and from $(\ref{SN48})$,
\bel{SN48X}|\bar v(t)-v(t,\gs)|\leq K_6\quad \text{for } (t,\gs)\in [T,\infty)\ti S^{N-1}.
\ee
The function $v$ satisfies $(\ref{SN7})$. Hence we apply \rlemma {int-2} with $\eta(t)=t$ and obtain that  $|v_t(t,\gs)|+|\nabla_gv(t,\gs)|\leq Ct$ in $[1,\infty)\ti S^{N-1}$. From this point we encounter three possibilities:\smallskip

\nind  (i) $\dsps\liminf_{t\to\infty}\bar v(t)>-\infty$. In that case $v$ remains bounded in $[T,\infty)\ti S^{N-1}$. Then \rprop{conv} applies, thus there exists a solution 
$\gw$ of $(\ref{Na-14})$ such that $v(t,.)\to\gw$ as $t\to\infty$. It follows that 
$$\lim_{r\to 0}r^2e^{u(r,.)}=e^{\gw}\quad\text{uniformly on }S^{N-1}.
$$  
This contradicts $(\ref{SN29})$.
\smallskip

\nind  (ii) $\dsps\liminf_{t\to\infty}\bar v(t)=-\infty$ and $\dsps\limsup_{t\to\infty}\bar v(t)>-\infty$. Then there exists an increasing sequence $\{t_n\}$ converging to $\infty$ such that 
\bel{SN49*}\left\{\BA{lll}
\bar v_t(t_n)=0\\
\bar v(t_{2n})\to-\infty\\
\bar v(t_{2n+1})>\bar v(t_{2n})\\
\bar v\text { is nondecreasing on }[t_{2n},†t_{2n+1}].
\EA\right.\ee
Averaging $(\ref{SN7})$ it is clear that the term
$me^{(q-2)t}\left((v_t+2)^2+|\nabla' v|\right)^{\frac q2}$ tends to zero exponentially and for some $\gd>0$, we obtain from $(\ref{SN7})$, $(\ref{SN8})$ and $(\ref{SN48X})$ the following differential  inequality,
\bel{SN50}
\bar v_{tt}-(N-2)\bar v_t-2(N-2-\gd)+ke^{\bar v}\geq 0\quad\text{for }\,t\geq T.
\ee
We multiply this inequality by $\bar v_t(t)$ for $t\in [t_{2n},t_{2n+1}]$ and integrate. Then 
$$k\left(e^{\bar v(t_{2n+1})}-e^{\bar v(t_{2n})}\right)\geq 2(N-2-\gd)(\bar v(t_{2n+1})-\bar v(t_{2n}))+(N-2)\int_{t_{2n}}^{t_{2n+1}}\bar v_t^2dt.
$$
If $\bar v(t_{2n+1})-\bar v(t_{2n})$ is unbounded this contradicts the fact that $e^v$ is bounded. If $\bar v(t_{2n+1})-\bar v(t_{2n})$ remains bounded we have by Rolle's theorem
$$k\left(e^{\bar v(t_{2n+1})}-e^{\bar v(t_{2n})}\right)=ke^{\gth_n \bar v(t_{2n+1})+(1-\gth_n)\bar v(t_{2n})}(\bar v(t_{2n+1})-\bar v(t_{2n}))\geq 2(N-2-\gd)(\bar v(t_{2n+1})-\bar v(t_{2n}))$$
for some $\gth_n\in (0,1)$. Hence $ke^{\gth_n \bar v(t_{2n+1})+(1-\gth_n)\bar v(t_{2n})}\geq 2(N-2-\gd)$ which is impossible since $\gth_n \bar v(t_{2n+1})+(1-\gth_n)\bar v(t_{2n})\to-\infty$.
As a consequence the assumption (ii) yields a contradiction again.\smallskip

\nind  (iii) $\dsps\lim_{t\to\infty}\bar v(t)=-\infty$. Then there exists a function $\lambda(t)$ tending to $0$ when $t\to\infty$ such that 
\bel{SN51}
\bar v_{tt}-(N-2)\bar v_t-2(N-2-\lambda(t))= 0\quad\text{for }\,t\geq T.
\ee
Since $\dsps\lim_{t_n\to\infty}e^{(2-N)t_n}\bar v_t(t_n)=0$ for some sequence $\{t_n\}$ tending to $\infty$, we obtain by integration from $t$ to $t_n$ and letting $\{t_n\}\to\infty$
\bel{SN51*}
-\bar v_t(t)=2e^{(N-2)t}\int_t^\infty\frac{N-2-\lambda(s)}{N-2}e^{(2-N)s}=2(1+o(1))\quad\text{as }t\to\infty,
\ee
which implies $\dsps\lim_{t\to\infty}\bar v_t(t)=-2$. Hence $\bar v(t)=-2t(1+o(1))$. Plugging this estimate into $(\ref{SN51})$ with the expression of the function $\lambda(t)$, it infers that
$|\lambda(t)|\leq ce^{-\gth t}$ for some $\gth>0$. Hence expression $(\ref{SN51*})$ can be made more precise as follows
$$
\bar v_t(t)=-2+O(e^{-\gth t}),
$$
which implies that 
$$\bar v(t)=-2t+K_2+O(e^{-\gth t})
$$
Hence if (iii) holds $u(x)$ remains bounded near $0$ which contradicts $(\ref{SN49})$.\smallskip

\nind As a consequence the claim of Step 2 holds. \smallskip

From Step 1-Step 2 the function $u$ is bounded and by standard regularity theory (e.g. \cite {GiTr}) it is a smooth solution in $B_1$.

\qeda

\subsection{The case $\frac{N}{N-1}\leq q<2$}


\blemma{conv-3} Assume that $\frac{N}{N-1}<q<2$ and $u$ is a solution of $(\ref{Na-1})$ in $B_1\setminus\{0\}$ bounded  from above near $0$. 
Then either $u$ can be extended as a $C^2$ solution of $(\ref{Na-1})$ in $B_1$, or  there exist positive constants $c_1,c_2$ such that
\bel{SN53}\dsps
-c_2|x|^{-\gb}\leq u(x)\leq -c_1|x|^{-\gb}\quad\text{for all }x\in B_1\setminus\{0\},
\ee
where $\gb=\frac{2-q}{q-1}$.
\es
\Proof Since $u$ is bounded from above we can apply \rcor{Harn1}: there exist $C>0$ and $a\in\BBR$ such that 
\bel{SN52}u(x)-a\leq C(u(y)-a)\quad\text{for all }x,y\text{ s.t. }0<|x|=|y|\leq\frac 12.
\ee
If $u$ is also bounded from below near $0$, it follows that $u$ is locally bounded and by classical regularity results 
it can be extended by continuity as a $C^2$ solution in $B_1$. If $u$ is not bounded from below the function $U=-u$ which satisfies 
 \bel{SN31}-\Gd U-m|\nabla U|^q+e^{-U}=0
\ee
and is bounded from below by assumption but not upper bounded. Assume now that 
\bel{SN54}\dsps
\liminf_{x\to 0}|x|^\gb U(x)=0.
\ee
It follows by $(\ref{SN52})$ that 
\bel{SN54*}\dsps\lim_{r_n\to 0}r_n^\gb\max_{|y|=r_n}U(y)=0,
\ee
for some sequence $\{r_n\}$ tending to $0$. Since $q>\frac{N}{N-1}$ the function
\bel{SN55}\dsps
x\mapsto U_\gb(x)=\Gl_{N,q}|x|^{-\gb}\quad\text{where }\Gl_{N,q}=\frac 1\gb\left(\frac{N-2-\gb}{N}\right)^{\frac 1\gb}=\frac 1\gb\left(\frac{(N-1)q-N}{q-1}\right)^{\frac 1{q-1}}
\ee
satisfies $-\Gd U_\gb-m|\nabla U_\gb|^q =0$ in $\BBR^N\setminus\{0\}$. Thus, for $\ge>0$ small enough, $U_{\gb,\ge}:=\ge U_\gb$ satisfies
\bel{SN55*}-\Gd U_{\gb,\ge}-m|\nabla U_{\gb,\ge}|^q=m\ge(1-\ge^{q-1})|\nabla U_{\gb,\ge}|^q\geq 0.
\ee
Since $\dsps \limsup_{x\to 0}	 U(x)=\infty$ as $x\to 0$, there exists $\gx\in B_1\setminus\{0\}$ such that $U(\gx)-\max_{|y|=1}U(y)> \frac 12$. 
There exists $\ge_0>0$ such that for $0<\ge<\ge_0$ there holds $U_{\gb,\ge}(\xi)< \frac 12$.
From $(\ref{SN54*})$, for any $0<\ge\leq \ge_0$, there exists $n_\ge\in\BBN$ such that for any $n\geq n_\ge$, we have $U(y)\leq U_{\gb,\ge}(y)$ for any $|y|=r_n$. The function 
$$x\mapsto W(x):=U(x)-U_{\gb,\ge}(x)-\max_{|y|=1}U(y)$$
satisfies
\bel{SN55**}-\Gd W-m\left(|\nabla U|^q-|\nabla U_{\gb,\ge}|^q\right)\leq -e^{-U}
\ee
in $B_1\setminus B_{r_n}$ and it is negative for $|x|=1$ and $|x|=r_n$. We can assume that $r_n<|\xi|$, hence exists 
 $W(\xi)>0$.
Therefore the maximum of $W$ in $B_1\setminus\overline {B_{r_n}}$
is achieved at some interior point $x_0$. Since $\nabla W(x_0)=0$ and thus $|\nabla U(x_0)|^q-|\nabla U_{\gb,\ge}(x_0)|^q=0$,  we have that
$$0\leq -\Gd W(x_0)\leq -e^{-U(x_0)}<0.
$$
This is a contradiction. Therefore $(\ref{SN54})$ cannot hold and we have
\bel{SN56}\dsps
U(x)\geq c_2|x|^{-\gb}.
\ee
This implies the claim.\qeda\medskip


The case $q=\frac N{N-1}$ necessitates an improvement of the estimate $(\ref{Na-19})$ obtained in \rcor{Th1}.

\blemma{conv-4} 
Assume that $q=\frac{N}{N-1}$ and $u$ is a solution of $(\ref{Na-1})$ in $B_1\setminus\{0\}$ upper bounded near $0$. Then either $u$ can be extended as a $C^2$ solution of $(\ref{Na-1})$ in $B_1$, or   there exist positive constants $c_3,c_4$ such that
\bel{SN58}\dsps
-c_4|x|^{2-N}\left(-\ln |x|\right)^{1-N}\leq u(x)\leq -c_3|x|^{2-N}\left(-\ln |x|\right)^{1-N}\quad\text{for all }x\in B_{\frac 12}\setminus\{0\}.
\ee
\es
\Proof As in the previous lemma we can assume that $u$ is not bounded from below, otherwise $u$ can be extended as a $C^2$ solution in $B_1$. Again the function $U=-u$ which satisfies $(\ref{SN31})$ is bounded from below but not from above. 
The equation 
\bel{SN59}
-\Gd U-m|\nabla U|^{\frac{N}{N-1}}=0
\ee
admits a radial singular solution $U^*$ which verifies
\bel{SN60}
U^*(r)=\frac{1}{N-2}\left(\frac{m}{N-1}\right)^{1-N}r^{2-N}(-\ln r)^{1-N}(1+o(1))\quad\text{as }r\to 0,
\ee
and
\bel{SN61}
U_r^*(r)=-\left(\frac{m}{N-1}\right)^{1-N}r^{1-N}(-\ln r)^{1-N}(1+o(1))\quad\text{as }r\to 0.
\ee
Furthermore we can assume that $U^*(1)=0$.\smallskip

\nind {\it Step 1.} We claim that
\bel{SN60-1}
U(x)\geq c_3U^*(|x|)\quad\text{if }0<|x|\leq \frac 12.
\ee
By contradiction, we assume that 
\bel{SN64}
\liminf_{x\to 0}\frac{U(x)}{U^*(|x|)}=0.
\ee
There exists a sequence $\{r_n\}$ tending to $0$ such that 
$$\lim_{n\to\infty}\frac{U(x)}{U^*(r_n)}=0\quad\text{uniformly for }|x|=r_n,
$$
the uniformity being still a consequence of Harnack inequality. For $\ge>0$ small enough, $U^*_{\ge}:=\ge U^*$ satisfies $(\ref{SN55*})$ with $q=\frac{N}{N-1}$. Again the function $\dsps W(x)=U(x)-U^*_{\ge}-\max_{|y|=1}U(y)$ satisfies $(\ref{SN55**})$ with $q=\frac{N}{N-1}$. Since $U$ is not upper bounded near $0$ there exists $\xi$ in $B_1\setminus\{0\}$ and $\ge_0>0$ such that  $W(\xi)>0$ for $0<\ge<\ge_0$. We conclude that it is impossible as in the proof of \rlemma{conv-3}. Hence we have $(\ref{SN60-1})$.
\smallskip

\nind {\it Step 2.} We claim that there exists $c_4>0$ such that
\bel{SN60-2}
U(x)\leq c_4U^*(|x|)\quad\text{if }0<|x|\leq \frac 12.
\ee
We proceed again by contradiction in assuming that
\bel{SN64-1}
\liminf_{x\to 0}\frac{U^*(|x|)}{U(x)}=0.
\ee
By the Harnack inequality there exists a sequence $\{r_n\}$ converging to $0$ such that 
\bel{SN64-2}
\liminf_{n\to \infty}\frac{U^*(r_n)}{U(x)}=0\quad\text{uniformly for }|x|=r_n.
\ee
From $(\ref{SN60})$ and $(\ref{SN61})$, for any $c>0$ there exists $R_{c}\in (0,1)$ such that 
$$m|\nabla U^*(|z|)| ^{\frac N{N-1}}\geq e^{-cU^*(|z|)}\quad\text{for all }z\in B_{R_{c}}\setminus \{0\}.
$$
For $\ge>0$, we set $\tilde  W(x):= U^*(|x|)-\ge U(x)-U^*(R_c)$. The function $\tilde W$ is negative on $\prt B_{R_c}$. Since  $U^*(|x|)\to\infty$ when $x\to 0$, there exists $\xi\in B_{R_c}\setminus \{0\}$ such that $U^*(|\xi|)-U^*(R_c)>\frac 12$, hence there exists $\ge_0>0$ such that $\ge U(\xi)<\frac 12$ for $\ge\leq\ge_0$ and thus $\tilde W(\xi)>0$. We fix $\ge\in (0,\ge_0)$. Because of $(\ref{SN64-2})$ one can find $n_0\in\BBN_*$ such that for $n\geq n_0$ one has
$r_n<|\xi|$ and $\tilde  W(x)<0$ for $|x|=r_n$. Therefore the maximum of $\tilde W$ in $B_{R_c}\setminus B_{r_n}$ is achieved at some interior point $x_0$. 
Since 
$$-\Gd(\ge U)-m|\nabla(\ge U)|^q=\ge\left(-\Gd U-m\ge^{q-1}|\nabla U|^q\right)=\ge m(1-\ge^{q-1})|\nabla U|^q-\ge e^{-U},
$$
and $e^{-U}\leq e^{-c_3U^*}$  from $(\ref{SN60-1})$, we obtain  
\bel{S65}\BA{lll}-\Gd \tilde W(x_0)-m\left(|\nabla U^*(|x_0|)|^{\frac N{N-1}}-|\nabla (\ge U(x_0))|^{\frac N{N-1}}\right)\\[2mm]
\phantom{----------------}=\ge(\ge^{{\frac 1{N-1}}}-1)|\nabla U(x_0)|^{\frac N{N-1}}+\ge e^{-U(x_0)}\\[2mm]
\phantom{----------------}
=m\ge^{-\frac 1{N-1}}(\ge^{\frac 1{N-1}}-1)|\nabla U^*(|x_0|)|^{\frac N{N-1}}+\ge e^{-U(x_0)}\\[2mm]
\phantom{----------------}
\leq m(1-\ge^{-\frac 1{N-1}})|\nabla U^*(|x_0|)|^{\frac N{N-1}}+\ge e^{-c_3U^*(|x_0|)}\\[2mm]
\phantom{----------------}
\leq m(1-\ge^{-\frac 1{N-1}}+\ge)|\nabla U^*(|x_0|)|^{\frac N{N-1}}.
\EA\ee
Up to changing $\ge_0$ we can assume that $1-\ge^{{-\frac 1{N-1}}}+\ge<0$ if  $\ge \in (0,\ge_0)$.  We obtain
\bel{SN66}
0\leq m(1-\ge^{-\frac 1{N-1}}+\ge)|\nabla U^*(|x_0|)|^{\frac N{N-1}}<0,
\ee
contradiction. Hence $(\ref{SN64})$ cannot hold, therefore $\frac{U^*(|x|)}{U(x)}\geq c_5$ and we obtain $(\ref{SN58})$.\qeda
\subsection{The radial case of \rth{Th5}}
We present here the complement of the proof of \rth{Th5}. 

\blemma{rad1} Let $1<q<2$ and $u(r)$ be a radial solution $(\ref{Na-1})$ in $B_1\setminus\{0\}$. Then $u_r(r)$ is monotone 
and we encounter the three following possibilities\smallskip

\nind (i) $u(r)\to\infty$ when $r\to 0$\smallskip

\nind (ii)  $u(r)\to \ell\in \BBR$ when $r\to 0$\smallskip

\nind (iii)  $u(r)\to-\infty$ when $r\to 0$.

\es
\Proof Since 
\bel{Rd0}
-u_{rr}-\frac{N-1}{r}u_r+m|u_r|^q-e^u=0
\ee
in $(0,1]$, $u$ cannot have any local minimum. Hence it is monotone, either increasing or decreasing and the claim follows.\qeda\\

\blemma{rad2} Let $\frac N{N-1}\leq q<2$ and $u(r)$ be a radial solution $(\ref{Na-1})$ in $B_1\setminus\{0\}$ such that $u(r)\to-\infty$ when $r\to 0$, then\smallskip

\nind(i) If $1<q<\frac {N}{N-1}$, relation $(\ref{Na-23})$ holds for some $\gg<0$.\smallskip

\nind(ii) If $\frac N{N-1}<q<2$, relation $(\ref{Na-27})$ holds.
\smallskip

\nind(iii) If $q=\frac N{N-1}$, relation $(\ref{Na-28})$ holds.
\es
\nind\Proof By \rth{Th4} the quantity $r^2e^{u(r)}$ is uniformly bounded in $B_1\setminus\{0\}$. Hence the result of \rth{Th5} applies. We give also below a simpler proof of \rth{Th5}-(iii) based upon ODE techniques.\\
 From monotonicity, we have that $e^{u(r)}\downarrow 0$  when $r\to 0$. Therefore $\dsps\limsup_{r\to 0} u_r(r)=\infty$. Set $h=u_r$, then 
$h\geq 0$ and 
\bel{Rd01}-h_{rr}-\frac{N-1}{r}h_{r}+\frac{N-1}{r^2}h+mq|u_{r}|^{q-1}h-e^{u}h=0.
\ee
If the function $h$ is not monotone near $0$, then there exists a sequence ${r_n}$ decreasing to $0$ such that $h(r_n)$ is a local maximum. Hence 
$$h(r_n)\left(\frac{N-1}{r_n^2}-e^{u(r_n)}\right)=h_{rr}(r_n)\leq 0.
$$ 
Because $\frac{N-1}{r_n^2}-e^{u(r_n)}\to\infty$ when $n\to\infty$, the only possibility is $h(r_n)=0$. Since $\eta_{r}(r_n)=0$ and $h$ satisfies a linear equation, it follows by the Cauchy-Lipschitz theorem that $u_r=h=0$ on $(0,1)$ and therefore $u$ is constant, contradiction. This implies that $u_r(r)$ is monotone and necessarily it satisfies
\bel{Rd1}
\lim_{r\to 0}u_r(r)=\infty.
\ee
Therefore, there exists a positive function $\tilde m:=\tilde m(r)$ satisfying $\tilde m(r)\to m$ when $r\to 0$, such that 
$$-u_{rr}-\frac{N-1}{r}u_r+\tilde mu_r^q=0\Longleftrightarrow -\frac{u_{rr}}{u_r^q}-\frac{N-1}{r}u_r^{1-q}+\tilde m=0.
$$
Set $X(r)=u_r^{1-q}$, then 
$$X'-\frac{(N-1)(q-1)}{r}X+(q-1)\tilde m=0\Longleftrightarrow\frac{d}{dr}\left(r^{-(N-1)(q-1)}X\right)+r^{-(N-1)(q-1)}(q-1)\tilde m=0.
$$
 It follows
\bel{Rd2}\BA{lll}\dsps r^{-(N-1)(q-1)}X(r)=r_0^{-(N-1)(q-1)}X(r_0)+(q-1)\int_r^{r_0}s^{-(N-1)(q-1)}\tilde m(s)ds.
\EA\ee
\smallskip

\nind {\it 1- Case $1<q<\frac{N}{N-1}$.} From the equation satisfied by $X$, we have
\bel{Rd2'}\BA{lll}\dsps
r^{-(N-1)(q-1)}X(r)+(q-1)\int_0^{r}s^{-(N-1)(q-1)}\tilde m(s)ds\\[4mm]
\phantom{------------}\dsps=r_0^{-(N-1)(q-1)}X(r_0)+(q-1)\int_0^{r_0}s^{-(N-1)(q-1)}\tilde m(s)ds.
\EA\ee
Therefore, if we set  
$$A_0:=r_0^{-(N-1)(q-1)}X(r_0)+(q-1)\int_0^{r_0}s^{-(N-1)(q-1)}\tilde m(s)ds\neq 0,
$$
we obtain
$$\BA{lll}\dsps
r^{-(N-1)(q-1)}X(r)=A_0+(q-1)\int_0^{r}s^{-(N-1)(q-1)}\tilde m(s)ds=A_0+o(1).
\EA$$
Since $X(r)=(u_r(r))^{1-q}$, it implies that 
\bel{Rd4}\BA{lll}\dsps
u_r(r)=A_0^{-\frac{1}{q-1}}r^{1-N}(1+o(1)).
\EA\ee
which in turn implies that $$u(x)=-\frac{A_0^{-\frac{1}{q-1}}}{N-2}r^{2-N}(1+o(1))$$
 as $r\to 0$.\smallskip

\nind {\it 2- Case $\frac{N}{N-1}\leq q<2$}. Then formula $(\ref{Rd2'})$ does not hold and it has to be replaced by $(\ref{Rd2})$ from which follows
$$X(r)=\left(\frac{r}{r_0}\right)^{(N-1)(q-1)}X(r_0)+(q-1)r^{(N-1)(q-1)}\int_r^{r_0}s^{-(N-1)(q-1)}\tilde m(s)ds,
$$
and by l'Hospital's rule
\bel{Rd5x}(q-1)r^{(N-1)(q-1)}\int_r^{r_0}s^{-(N-1)(q-1)}\tilde m(s)ds=\left\{\BA{lll}\dsps\frac{(q-1)mr}{(N-1)(q-1)-1}(1+o(1))\quad&\text{if }q>\frac N{N-1}\\[4mm]\dsps
(q-1)mr\ln\frac{r_0}{r}(1+o(1))\quad&\text{if }q=\frac N{N-1}
\EA\right.
\ee
as $r\to 0$. Hence 
\bel{Rd5}\BA{lll}\dsps
u_r(r)=\left\{\BA{lll}\dsps\left(\frac{(N-1)(q-1)-1}{(q-1)mr}\right)^{\frac {1}{q-1}}(1+o(1))\quad&\text{if }q>\frac N{N-1}\\[4mm]\dsps
\left((q-1)mr\ln\frac{r_0}{r}\right)^{-\frac 1{q-1} }(1+o(1))\quad&\text{if }q=\frac N{N-1}.\EA\right.
\EA\ee
Thus we obtain the following asymptotics as $r\to 0$,
\bel{Rd6}\BA{lll}\dsps
u(r)=\left\{\BA{lll}\dsps-\left(\frac{(N-1)(q-1)-1}{(q-1)m}\right)^{\frac {1}{q-1}}\frac {r^{-\gb}}{\gb}(1+o(1))\quad&\text{if }q>\frac N{N-1}\\[4mm]\dsps
\frac{(q-1)^{-\gb}}{2-q}r^{-\gb}\left(\ln \frac{r_0}{r}\right)^{-\gb-1}(1+o(1))\quad&\text{if }q=\frac N{N-1},\EA\right.
\EA\ee
where $\beta$ is the exponent defined at $(\ref{Na-6})$. This ends the proof.
\qeda 


\mysection{The case $q>2$}

We have the following form of Harnack inequality, which is actually not optimal since the a priori estimate on $u$ is not the general one which follows from \rth{Th3}.

\bprop{q2-1} Assume $N\geq 2$ and $q>2$. Let $u$ be a solution of $(\ref{Na-1})$ in $B_1\setminus\{0\}$ satisfying 
\bel{Q2-1}
\lim_{x\to 0}|x|^{q}e^{u}=0.
\ee
Then $u$ can be extended as a continuous function in   $B_1$ which satisfies
\bel{Q2-1*}
|u(x)-u(0)|\leq C|x|^{\frac{q-2}{q-1}}.
\ee
Furthermore, if $u$ is radial, then either $u_r(0)=0$ and $u$ is a regular solution, or 
\bel{Q2-1**}
u_r(r)=\left(\frac{N(q-1)-q}{m(q-1)r}\right)^{\frac{1}{q-1}}(1+o(1))\quad\text{as }r\to 0.
\ee
\es
\Proof 
{\it Step 1.} We first prove that the result holds if 
\bel{Q2-2}
|x|^{\gk}e^{u}\leq K\quad\text{in }B_1\setminus\{0\},
\ee
for some $0<\gk<2$ and $K>0$.
Indeed from \rth{Th2}, we have
$$|\nabla u(x)|\leq C\left(|x|^{-\frac{1}{q-1}}+|x|^{-\frac\gk q}\right).
$$
Therefore 
$$m|\nabla u(x)|^q+e^u\leq C\left(|x|^{-\frac q{q-1}}+|x|^{-\gk}\right).
$$
Since $q>2$ and  $\gk<2$, $|x|^{-\frac q{q-1}}+|x|^{-\gk}\in L^p(B_1)$ for some $p>\frac N2$. Hence $u\in W^{2,p}(B_\frac 12)\subset C^{0,\gth}(B_\frac 12)$ for some $\gth\in (0,1)$.\smallskip

\nind {\it Step 2. H\"older continuity.} Under assumption $(\ref{Q2-1})$, for any $\ge>0$ there exist $C_\ge>0$ such that 
$$e^{u(x)}\leq \ge|x|^{-q}+C_\ge\quad\text{for }0<|x|<1.
$$
Using again  \rth{Th2} we obtain
$$|\nabla u(x)|\leq C\left(|x|^{-\frac{1}{q-1}}+\ge^\frac{1}{q} |x|^{-1}+C_\ge^{\frac1q}\right) \quad\text{for }0<|x|<\frac12,
$$
and by integration
$$|u(x)|\leq \max_{|z|=\frac12}|u(z)|+\ge^\frac{1}{q} C\ln\frac{1}{2|x|}+C'.
$$
This implies that 
$$|x|^{\ge^\frac{1}{q} C}e^{u(x)}\leq K'
$$
for some $K'>0$. For $\ge$ small enough we have $C\ge^\frac{1}{q}<2$ and we conclude with Step 1.$\phantom{------------}$.\smallskip

\nind {\it Step 3. Proof of $(\ref{Q2-1*})$.} Since $e^u$ is bounded by Step 1, we have by \rth{Th2}
$$|\nabla u(x)|\leq C(|x|^{-\frac 1{q-1}}+1).
$$
By integration we obtain $(\ref{Q2-1*})$.
\smallskip

\nind {\it Step 4. The radial case.} If $u$ is radial we have
$$-u_{rr}-\frac{N-1}{r}u_r+m|u_r|^q-e^{u}=0.
$$
Then either $u$ is a regular solution and hence $u_r(0)=0$, or $u$ is singular in the sense that $u_r$ is not bounded near $0$. At each $r_0$ such that $u_r(r_0)=0$, we have $u_{rr}(r_0)<0$, hence $u$ cannot oscillate since $u_r$ keeps a constant sign near $0$. By the argument of the proof of \rlemma {rad2} the function $u_r$ is monotone near $0$. Since $u_r$ is not bounded, it implies that 
$|u_r(r)|\to\infty$. Hence for $\ge>0$ we have 
$$\dsps
-u_{rr}-\frac{N-1}{r}u_r+(m+\ge)|u_r|^q\leq 0,
$$
and 
$$
\dsps-u_{rr}-\frac{N-1}{r}u_r+(m-\ge)|u_r|^q\geq 0.
$$
Integrating the two inequalities we have that $u_r$ is positive and furthermore
\bel{Q2-2*}
\left(\frac{N(q-1)-q}{(q-1)(m+\ge)r}\right)^{\frac{1}{q-1}}
\leq u_r(r)\leq \left(\frac{N(q-1)-q}{(q-1)(m-\ge)r}\right)^{\frac{1}{q-1}}.
\ee
This implies $(\ref{Q2-1**})$ since $\ge$ is arbitrary.
$\phantom{------------}$\qeda\medskip

\nind\Remark (i) The condition $(\ref{Q2-1})$ can be slightly improved and replace by 
\bel{Q2-2***}
\limsup_{x\to 0}|x|^{q}e^{u}=\ge_0<\left(\frac 2C\right)^q,
\ee
where $C$ is the constant appearing in the gradient estimate from \rth{Th2}.\smallskip

\nind (ii) We conjecture that the radiallity can be dropped in the proof of $(\ref{Q2-1**})$ and that the strongest estimate holds for singular solutions
\bel{Q2-2****}
|u(x)-u(0)|=\frac{q-1}{q-2}\left(\frac{N(q-1)-q}{(q-1)m}\right)^{\frac{1}{q-1}}|x|^{\frac{q-2}{q-1}}(1+o(1))\quad\text{as }x\to 0.
\ee

\bth{clok}  Assume $N\geq 2$ and $q>2$. If $u$ is a solution of $(\ref{Na-1})$ in $B_{1}\setminus\{0\}$.
satisfying
\bel{Q2-2n}\dsps \liminf_{x\to 0}|x|^qe^{u(x)}>0,\ee
then there holds
\bel{Q2-15}
u(x)=q\ln\frac {1}{|x|}+\ln mq^q+o(1)\quad\text{as }x\to 0.
\ee
Furthermore, if $u$ is a radial function,  relation $(\ref{Q2-15})$ still holds if we assume $\dsps\limsup_{x\to 0}r^qe^{u(r)}>0$ instead of inequality $(\ref{Q2-2n})$.
\es
\Proof {\it Step 1: The general case}: Assume that $u$ is non-necessary radial and satisfies
\bel{Q2-16}
\liminf_{x\to 0}|x|^qe^{u(x)}=\gk>0.
\ee
For $\ell>0$ the function  $T^{ei}_\ell[u]$ satisfies for some $,\gk_2,\gk_1>0$,
$$q\ln\frac 1{|x|}+\ln\gk_2\geq T^{ei}_\ell[u](x)=u(\ell x)+q\ln\ell\geq q\ln\frac 1{|x|}+\ln\gk_1\quad\text{for all }x\in B_{\ell^{-1}},
$$
an inequality which follows from $(\ref{Q2-16})$. Furthermore \rth{Th3} implies  also
$$
|x|^qe^{T_\ell^{ei}[u](x)}+|x||\nabla T_\ell^{ei}[u](x)|\leq C\quad\text{for all }x\in B_{(2\ell)^{-1}}\setminus\{0\},
$$
which yields
$$T^{ei}_\ell[u](x)\leq q\ln\frac 1{|x|}+\ln C,$$
 and we have that
$$\ell^{q-2}\Gd T_\ell^{ei}[u]+m|\nabla T_\ell^{ei}[u]|^q-e^{T_\ell^{ei}[u]}=0\quad\text{in }B_{\ell^{-1}}\setminus\{0\}.
$$
We apply the Crandall-Lions result (see \cite{CrLi}) in the domain $B_n\setminus B_{\frac 1n}$: for $\ell <\frac 1{2n}$, the set of functions $\{T^{ei}_\ell[u]\}$
is relatively compact in $C^{0,\gth}(\overline {B_n\setminus B_{\frac 1n}})$ for some $\gth\in (0,1)$. Hence there exist a sequence $\{\ell_{k,n}\}_{k\in\BBN}$ tending to $0$ when $k\to\infty$ with 
$\ell_{k,n} <\frac 1{2n}$ and a function $\hat u_n\in C^{0,\gth}(\overline {B_n\setminus B_{\frac 1n}})$ which is a viscosity solution of 
$$m|\nabla \hat u_n|^q-e^{ \hat u_n}=0\quad\text{in }B_n\setminus B_{\frac 1n},
$$
such that $\{T_{\ell_{k,n}}^{ei}[u]\}$  converges  to $\hat u_n$ in the $C^{0,\gth}$ topology of $\overline {B_n\setminus B_{\frac 1n}})$. By Cantor diagonal sequence there exists 
a sequence $\{\ell_{\phi(n),n}\}_{n\in\BBN}$ where $\phi(.)$ is an increasing entire function and a function $\hat u$ such that $\{T_{\ell_{\phi(n),n}}^{ei}[u]\}$  converges  to $\hat u$
in the $C^{0,\gth}_{loc}$ topology of $\BBR^N\setminus\{0\}$, and $\hat u$ is a  viscosity solution $\hat u$ of 
$$m|\nabla \hat u|^q-e^{\hat u}=0\quad\text{in }\BBR^N\setminus\{0\},
$$
such that $T_{\ell_n}^{ei}[u]\to \hat u$ locally uniformly in $\BBR^N\setminus\{0\}$ and $\nabla T_{\ell_n}^{ei}[u]$ to $\nabla \hat u$ in the local weak-star topology of $L^\infty(\BBR^N\setminus\{0\})$. Furthermore there holds in $\BBR^N\setminus\{0\}$,
$$
|x|^qe^{\hat u(x)}+|x||\nabla \hat u(x)|\leq C,
$$
and
$$q\ln\frac 1{|x|}+\ln\gk_2\geq \hat u(x)\geq q\ln\frac 1{|x|}+\ln\gk_1.
$$
We set $v(x)=e^{-\frac{\hat u(x)}{q}}$. Then 
$$|\nabla v(x)|=\frac{1}{qm^\frac 1q}\quad\text{in }\BBR^N\setminus\{0\}.
$$
Hence $v$ is Lipschitz continuous and it is classical (see e.g. \cite{Ca-Cr}) that
\bel{Q2-16*}
v(x)-v(0)=\frac{|x|}{qm^\frac 1q}\Longrightarrow \hat u(x)=q\ln\frac{qm^{\frac 1q}}{|x|},
\ee
here we use the fact that that $v(0)=0$. The uniqueness of $\hat u$ implies that $T_\ell^{ei}[u]\to \hat u$ locally uniformly $\BBR^N\setminus\{0\}$ as $\ell\to 0$,
which reads
$$u(\ell x)=q\ln\frac{1}{\ell^q|x|^q}+\ln mq^q+o(1).
$$
Taking $|x|=1$ and $\ell=|y|$ yields

\bel{Q2-16**}
u(x)=\ln\frac{1}{|x|^q}+\ln mq^q+o(1)\quad\text{as }x\to 0.
\ee

\nind {\it Step 2: The radial case}: let $u$ be a radial solution such that 
\bel{Q2-17}
\limsup_{r\to 0}r^qe^{u(r)}=\gk>0.
\ee
By \rlemma{erra} $u$ is monotone near $r=0$ and actually decreasing because of $(\ref{Q2-17})$. For $c\in (0,1)$ we set $F_c(r)=e^{u(r)}-mc|u_r|^q$. If $F_c$ does not keep a constant sign, near $0$, 
there exists $\{r_n\}$ tending to $0$ such that $F_c(r_n)=0$ and $F'_c(r_n)\leq 0$. Since $F'_c(r)=u_r\left(e^{u}-qmc|u_r|^{q-2}u_{rr}\right)$, we have $e^{u(r_n)}=mc|u_r(r_n)|^q$, hence
\bel{Q2-17x}F'_c(r_n)=-mc|u_r(r_n)|^{q-1}\left(|u_r(r_n)|-\frac{q(N-1)}{r_n}-qm(1-c)|u_r(r_n)|^{q-1}\right).
\ee
Since $r_n\to 0$, we have that $F'_c(r_n)>0$, contradiction. Hence $F_c$ keeps a constant sign near $0$. By \rth{Th3} we have that 
$r^q|u_r(r)|^{q}\leq C^q$, therefore
$$r^qF_c(r)\geq r^qe^u-mcr^q|u_r(r)|^{q}\geq r^qe^u-cmC^q.
$$
If $cmC^q\leq \frac\gk2$ it follows that there exists a sequence $\{r'_n\}$ tending to $0$ such that $F_c(r'_n)>0$ and thus $F_c(r)>0$ near $0$. This yields that 
differential inequality
$$e^{u}\geq mc|u_r|^{q}\Longrightarrow \left(e^{-\frac uq}\right)_r\leq\frac{1}{q(mc)\frac1q}\Longrightarrow u(r)\geq q\ln\frac{q(mc)^\frac1q}{r}.
$$
This implies $\dsps \liminf_{r\to 0}r^qe^{u(r)}>0$ and we conclude by Step 1.

\qeda\medskip

\nind\Remark As a consequence of $(\ref{Q2-15})$, we have
\bel{Q2-1xx}
\lim_{x\to 0}|x|^qe^{u(x)}=mq^q.
\ee\\

\nind\Remark The existence of singular solutions satisfying $(\ref{Q2-1*})$ is not difficult to obtain by a variant of the method of super and sub solutions. For the existence of solutions of eikonal type (actually radial), the proof of existence is much more difficult and is proved in \cite{BV-V3}.
\\

\nind{\it Open problem}. We conjecture that  the assumption$( \ref{Q2-2n})$ can be weaken and replaced by 
\bel{Q2-2n*}\dsps \limsup_{x\to 0}|x|^qe^{u(x)}>0.\ee
If it were true this would give a complete classification of singularities of solutions of $(\ref{Na-1})$ in $B_{1}\setminus\{0\}$ in the case $q>2$. Note that this result valid if $u$ is a radial function.

\mysection{Behaviour at infinity}
In this section we give upper estimates of solutions of $(\ref{Na-1})$ in an exterior domain of $\BBR^N$.
\subsection{Estimates of supersolutions}
We recall that if $u$ is a continuous function defined in $B_{r_0}^c$ we have denoted by $\gm(r)$ the infimum of $u(x)$ on $\{x:|x|=r\}$.
\bth{th*1*} Let $N\geq 2$, $q>1$ and $u$ be a supersolution of $(\ref{Na-1})$ in $B_{r_0}^c$. Then there exists $C>0$ depending on 
$N,m,q$ and $u$ such that 
\bel{Q4-1}
e^{\gm(r)}\leq Cr^{-\min\{2,q\}}\quad\text{for all }r\geq 2r_0.
\ee
\es
\Proof The proof is an adaptation of the one of \rth{Th4} to this framework, and we keep the same notations:
$ w=e^u$ and $\dsps M(r)=\min_{|x|=r}w(x)=e^{\gm(r)}$. 
 From \cite[Lemma4.4]{BV-V2} we know that the function $\gm$ is monotone for $r\geq \gr>r_0$, and so is the function $M$. We fix $R>\gr+1$. If $M$ is nondecreasing 
 on $[R-\ge,R+\ge]$, then $M(R)\geq w(\tilde x_{R,\ge})\geq M(|\tilde x_{R,\ge}|)\geq M((1-\ge)R)$, where we have denoted by $\tilde x_{R,\ge}$ the point of minimum of $w$ in $\overline B_{R+\ge}\setminus B_{R-\ge}$. Then
 $$M((1-\ge)R)\leq w(\tilde x_{R,\ge})\leq C\left(\frac{M^q(R)}{R^2}+\frac{M^q(R)}{R^q}\right)\leq\left\{\BA{lll}
 2C\frac{M^{\frac q{q+1}}(R)}{R^{\frac{2}{q+1}}}\;\text{ if }q>2\\[2mm]
  2C\frac{M^{\frac q{q+1}}(R)}{R^{\frac{q}{q+1}}}\;\text{ if }q<2.
 \EA\right.
 $$
 From the bootstrap Lemma
 \bel{Q4-2}
M(R)\leq \left\{\BA{lll}
C_1R^{-2}\;\text{ if }q>2\\[1mm]
2CR^{-q}\;\text{ if }q<2,
\EA\right.
\ee
which is in contradiction with the fact that $M$ is nondecreasing. Hence $M$ is nonincreasing. Therefore we have that
 $$M((1+\ge)R)\leq\left\{\BA{lll}
 2C\frac{M^{\frac q{q+1}}(R)}{R^{\frac{2}{q+1}}}\;\text{ if }q>2\\[2mm]
  2C\frac{M^{\frac q{q+1}}(R)}{R^{\frac{q}{q+1}}}\;\text{ if }q<2.
 \EA\right. $$
 Using again the bootstrap Lemma we derive  $(\ref{Q4-2})$ as before and $(\ref{Q4-1})$ finally.\qeda\medskip
 
 Since the spherical minimum of a global solution $u(x)$ tends to $-\infty$ when $|x|\to\infty$, we have the following Liouville type result.
 \bcor{Liou} If $N\geq 2$ and $q>1$ there exists no positive solution of $(\ref{Na-1})$ in $\BBR^N$.
 \es
\subsection{The case $1<q<2$}
The first result is an priori estimate which is the equivalent (both in its claim as in its proof) of \rth{Th3} valid for the solutions of $(\ref{Na-1})$ in an exterior domain which tend to $-\infty$.
\bth{th*1} Let $N\geq 2$, $1<q<2$ and $u$ be a solution of $(\ref{Na-1})$ in $B_1^c$ such that 
\bel{Q5-1}
\lim_{|x|\to \infty}e^{u(x)}=0.
\ee
Then there exists $C>0$ depending on $N$, $m$, $q$, and $u$ such that 
\bel{Q5-2}
e^{u(x)}\leq C|x|^{-q}\quad\text{and }\;|\nabla u(x)|\leq C|x|^{-1}\quad\text{in }B_2^c.
\ee
\es
\Proof We set $w=e^u$ and keep the notations of \rth{Th3}, assuming that $w\leq 1$. If we can prove that 
\bel{Q5-3}
w^\frac 1q(x)\leq \frac{C}{|x|},
\ee
then applying  $(\ref{Na-2-1})$ in $B_{\frac {|x|}{2}}(x)$ we have by \rth{Th2}
\bel{Q5-4}
|\nabla u(x)|\leq c_1|x|^{-\frac{1}{q-1}}+c_2\max_{z\in B_{\frac {|x|}{2}}(x)}w^{\frac 1q}+c_3\max_{z\in B_{\frac {|x|}{2}}(x)}w^{\frac {1}{2(q-1)}}.
\ee
Note that this estimate does not depend on the value of $q$ with respect to $2$. Since $q<2$, $|x|>1$ and $0<w\leq 1$, there holds $|x|^{-\frac{1}{q-1}}\leq |x|^{-1}$ and $w^{\frac {1}{2(q-1)}}\leq w^{\frac 1q}$. Hence $(\ref{Q5-3})$ implies 
$(\ref{Q5-2})$. Following the proof of \rth{Th3} we set 
$$M(x)=w^{\frac 1q}(x).
$$
By assumption $M(x)\to 0$ when $|x|\to\infty$. We apply again the doubling lemma with $X=\bar B_1$, $D=\BBR^N\setminus \bar B_2$, $\Gs=\BBR^N\setminus  B_2$, 
$\Gg=\prt B_2$ and $k=n$. There exists a sequence $\{y_n\}\subset \BBR^N\setminus \bar B_2$ such that $(|y_n|-2)M(y_n)\to\infty$ when $n\to\infty$.
Let $x_n\in \BBR^N\setminus \bar B_2$ such that $(\ref{Na-2*-5'''})$ holds. It is clear that $\{x_n\}$ is unbounded since $M$ is bounded on bounded subsets of 
$\BBR^N\setminus \bar B_2$, and we can assume that $|x_n|\to\infty$ when $n\to\infty$. Performing the same proof as in \rth{Th3} we obtain that 
$\ge_n=M^{2-q}(x_n)\to 0$ since $M(x_n)\to 0$ and $1<q<2$. The end of the proof is similar.\qeda\medskip

The following result is the counterpart of \rth{clok} at infinity.
\bth{clok-2}  Assume $N\geq 2$ and $1<q<2$. If $u$ is a solution of $(\ref{Na-1})$ in $B^c_{1}$
satisfying $(\ref{Q5-1})$ and 
\bel{Q2-2n**}\dsps \liminf_{|x|\to \infty}|x|^qe^{u(x)}>0.\ee
Then there holds
\bel{Q2-15*}
u(x)=q\ln\frac {1}{|x|}+\ln mq^q+o(1)\quad\text{as }|x|\to \infty.
\ee
Furthermore, if $u$ is radial relation $(\ref{Q2-15*})$ holds if we assume $\dsps\limsup_{|x|\to \infty}r^qe^{u(r)}>0$ instead of $(\ref{Q2-2n**})$.
\es
\Proof It follows from $(\ref{Q5-2})$ and $(\ref{Q2-2n**})$ that for $\ell>1$, the function $T^{ei}_\ell[u]$  satisfies
$$q\ln\frac1{|x|}+\ln\gk_1\leq T^{ei}_\ell[u](x)\leq q\ln\frac1{|x|}+\ln\gk_2,
$$
and
$$|x|^qe^{T^{ei}_\ell[u](x)}+|x||\nabla T^{ei}_\ell[u](x)|\leq C,
$$
in $B^c_{(2\ell)^{-1}}$, and that the following equation holds therein
$$\ell^{q-2}\Gd T^{ei}_\ell[u]+m|T^{ei}_\ell[u]|^q-e^{T^{ei}_\ell[u]}=0.
$$
 The conditions for applying Crandall-Lions result as in the proof of \rth{clok} are fulfilled and for $n<\frac\ell2$ there exists a sequence $\{\ell_{k,n}\}$ tending to $\infty$ when $k\to\infty$ and 
 a function $\hat u_n$ such that $\{T^{ei}_{\ell_{k,n}}\}$ converges in the $C^{0,\gth}(\overline{B_n\setminus B_{\frac 1n}})$ to $\hat u_n$ and $\hat u_n$ is the solution of the eikonal equation 
 $(\ref{Na-3})$. Using again a diagonal sequence $\{\ell_{\phi(n),n}\}$ and a function $\hat u$ such that $\{T^{ei}_{\ell_{\phi(n),n}}\}$ converges to $\hat u$ in the $C_{loc}^{0,\gth}(\BBR^N\setminus \{0\})$-topology and $\hat u$ satisfies  $(\ref{Na-3})$ in $\BBR^N\setminus \{0\}$ in which domain the following estimates hold
\bel{Qx-1}\BA{lll}\dsps
(i)\qquad\qquad  &|x|^qe^{\hat u(x)}+|x||\nabla \hat u(x)|\leq C\qquad\qquad\qquad\qquad\qquad\qquad\qquad\qquad\qquad\qquad\qquad\qquad\\[2mm]
 \dsps 
(ii) &q\ln\frac1{|x|}+\ln\gk_1\leq \hat u(x)\leq q\ln\frac1{|x|}+\ln\gk_2.
 \EA
\ee
Setting again $v(x)=e^{-\frac{\hat u(x)}q}$ and using the fact that $\hat u(0)=\infty$ from $(\ref{Qx-1})$, we conclude again that $(\ref{Q2-15*})$ holds and the conclusion follows as in \rth{clok}.\smallskip

\nind In the radial case, we know from \rlemma{erra} that $u(r)$ is monotone for $r\geq \gr_1>1$ and it is decreasing because of $(\ref{Q5-1})$. We define the function $F_c$ as in the proof of \rth{clok}. Since $u_r(r)\to 0$ when $r\to\infty$ and $1<q<2$, if $F_c$ does not keep a constant sign at infinity, and there exists a sequence $\{r_n\}$ tending to infinity and such that 
$F_c(r_n)=0$ and $F'_c(r_n)\leq 0$, the expression $(\ref{Q2-17x})$ shows that $F'_c(r_n)> 0$ which is a contradiction. Hence $F_c$ has a constant sign and for $c>0$ small enough we see that 
its sign is positive. Integrating the corresponding inequality between $1$ and $r>1$ we obtain
$$u(r)\geq q\ln\left(\frac{1}{q(mc)^\frac1{q}(r-1)+e^{-\frac{u(1)}{q}}}\right).
$$
This implies $(\ref{Q2-2n**})$.\qeda


\subsection{The case $q>2$}
In this section we prove the counterpart  for negative solutions of $(\ref{Na-1})$ in an exterior domain of the Harnack inequality and the isotropy property.
\bprop{th*2} Let $N\geq 2$, $q>2$ and $u$ be a solution of $(\ref{Na-1})$ in $B_1^c$ such that $|x|^{2}e^{u(x)}$ remains bounded in $B_R^c$,
Then there exists $0<C<1$ depending on $N, p, m$ and $\dsps\sup_{|x|\geq R}|x|^{2}e^{u(x)}$ such that 
\bel{Q6-2}
\max_{|x|=r}u(x)\leq C\min_{|x|=r}u(x)\quad\text{for all }r\geq R+1.
\ee
\es
\Proof We can assume that $u<0$ in $B_R^c$. We write $(\ref{Na-1})$ under the form
$$-\Gd u+A(x)u+B(x)|\nabla u|=0,
$$
where $A(x)=\frac{e^u}{-u}$ and $B(x)=m|\nabla u|^{q-1}$.
Assuming that $u\leq -1$, we have, 
$$|x|^2|A(x)|= |x|^2\frac{e^u}{|u|}\leq |x|^2e^{u(x)}\leq M.
$$
Next
$$\BA{lll}|x|B(x)\leq C\left(1+|x|e^{\frac {q-1}{q}\max\{u(z):|z-x|\leq\frac{|x|}{2}\}}+|x|e^{\frac {1}{2}\max\{u(z):|z-x|\leq\frac{|x|}{2}\}}\right).
\EA$$
If $q>2$ we have 
$$|x|e^{\frac {q-1}{q}\max\{u(z):|z-x|\leq\frac{|x|}{2}\}}+|x|e^{\frac {1}{2}\max\{u(z):|z-x|\leq\frac{|x|}{2}\}}\leq M'(1+|x|^{\frac{2-q}{q}})\leq 2M'.
$$
Therefore by \cite{GiTr}, Harnack inequality holds and $(\ref{Q6-2})$ follows.\qeda\medskip

The following result will be useful in the proof of the isotropy estimate

\blemma{th*2*} Let $N\geq 2$, $q>2$ and $u$ be a solution of $(\ref{Na-1})$ in $B_1^c$ such that $|x|^2e^{u(x)}$ remains bounded in $B_1^c$
Then there exists $C>0$ depending on $N$, $m$, $q$, and $M$ such that for all $x\in B_2^c$ there holds,
\bel{Q5-6}\BA{lll}
(i)\qquad\qquad\qquad\qquad&-C|x|^{\frac{q-2}{q-1}}\leq u(x)\leq -2\ln |x|+K_1\qquad\qquad\qquad\qquad\qquad\qquad\qquad\\[1mm]
(ii)\qquad\qquad&|\nabla u(x)|\leq C|x|^{-\frac{1}{q-1}}\quad\text{in }B_2^c.
\EA\ee
\es
\Proof The upper estimate (i) is the assumption. We write $(\ref{Q5-4})$ under the form
\bel{Q5-4^*}
|\nabla u(x)|\leq c_1|x|^{-\frac{1}{q-1}}+c_2\max_{z\in B_{\frac {|x|}{2}}(x)}e^{\frac {u(z)}q}+c_3\max_{z\in B_{\frac {|x|}{2}}(x)}e^{\frac {u(z)}{2(q-1)}}.
\ee
Since $q>2$ and $u(x)\to-\infty$ when $|x|\to\infty$, we have, assuming $u(x)\leq -1$, 
$$e^{\frac {u(x)}q}\leq e^{\frac {u(x)}{2(q-1)}} \leq M^{\frac1{2(q-1)}}|x|^{-\frac{1}{q-1}}\quad\text{for  }|x|\geq 2.
$$
Thus $(\ref{Q5-6})$-(ii) holds and by integration we obtain $(\ref{Q5-6})$-(i).\qeda\medskip

The isotropy estimate of $u-\bar u$ is obtained by a modified representation formula. We  set 
\bel{Q6-3*} t=\ln r\Longleftrightarrow r=e^t \quad\text{and }u(r,\gs)=v(t,\gs)-2\ln r\Longrightarrow v(t,\gs)=u(r,\gs)+2t.
\ee
Then $v$ satisfies in $[T,\infty)\ti S^{N-1}$
\bel{Q6-3}
v_{tt}+(N-2)v_t+\Gd_g v-2(N-2)-me^{(2-q)t}\left((v_t+2)^2+|\nabla_gv|^2\right)^{\frac q2}+e^v=0.
\ee
\bprop{th*3} Let $N\geq 2$, $q>2$ and $u$ be a negative solution of $(\ref{Na-1})$ in $B_1^c$ such that $|x|^2e^u$ is bounded.
 Then there exists $C_1>0$ depending on $N, p, m$ such that 
\bel{Q6-4}
\norm{u(r,.)-\bar u(r)}_{L^\infty(S^{N-1})}\leq C_1\quad\text{for all }r\geq 2.
\ee
\es
\Proof Since $|x|^{2}e^{u(x)}$ is bounded in $B_1^c$, estimate $(\ref{Q5-6})$  implies
\bel{Q6-4**}
\left((v_t+2)^2+|\nabla' v|^2\right)^\frac{1}{2}\leq Ce^{\gl t}:= Ce^{\frac{q-2}{q-1}t}\quad\text{in }[T,\infty)\ti S^{N-1}.
\ee
We define $\psi$ by $v(t,\gs)=e^{\frac{q-2}{q-1} t}\psi (t,\gs)=e^{\gl t}\psi (t,\gs)$, thus 
\bel{Q6-5}\BA{lll}
\psi_{tt}+(N-2+2\gl)\psi_t+\gl(N-2+\gl)\psi+\Gd_g\psi-2(N-2)e^{-\gl t}+e^{-\gl t}e^{e^{\gl t}\psi}\\[3mm]
\phantom{----------------}
-m\left((\psi_t+\gl\psi+2e^{-\gl t})^2+|\nabla_g\psi|^2\right)^{\frac q2}=0.
\EA\ee
If $\bar\psi(t)$ is the spherical average of $\psi(t,.)$, then $\psi^*=\psi-\bar \psi$ satisfies in $\BBR^+\ti S^{N-1}$
\bel{Q6-5+1}\BA{lll}
\psi^*_{tt}+a'\psi^*_t-\ell'\psi^*+\Gd_g\psi^*+F^*=0.
\EA\ee
where
\bel{Q6-5+2}\BA{lll}
a'=N-2+2\gl\quad\text{and }\,\ell'=\gl(2-N-\gl),
\EA\ee
and
where 
$$\BA{lll}F^*(t)=m\left(\left((\psi_t+\gl\psi+2e^{-\gl t})^2+|\nabla_g\psi|^2\right)^{\frac q2}-\overline{\left((\psi_t+\gl\psi+2e^{-\gl t})^2+|\nabla_g\psi|^2\right)^{\frac q2}}\right)\\[2mm]
\phantom{----------------------------}
-e^{-\gl t}\left(e^{e^{\gl t}\psi}-\overline{e^{e^{\gl t}\psi}}\right).
\EA$$
Noticing that $\frac{a'^2}{4}+\ell'=\frac{a^2}{4}+\ell=\frac{(N-2)^2}{4}$, we can 
 apply {\it formally} the representation formula used in the proof of \rprop{int-3} where the operator $\BBL$ and thus the semigroup $S(t)=e^{t\BBL}$ are unchanged.{\it Mutatis mutandis} we finally obtain that $\psi^*=\psi-\bar \psi$ is expressed by
\bel{Q6-6}\BA{lll}\dsps
\psi^*(t)=e^{(\frac{2-N}{2}-\gl)t}S(t)[\psi^*(0)]-\int_0^te^{(\frac{2-N}{2}-\gl)s}S(s)\left[\int_0^\infty e^{(\frac{N-2}{2}+\gl)\gt}S(\gt)[F^*(t+\gt-s)]d\gt \right]ds.
\EA\ee
In order to derive an upper estimate concerning the term $F^*(t+\gt-s))$ we need an upper estimate of the term $\left((\psi_t+\gl\psi+2e^{-\gl t})^2+|\nabla_g\psi|^2\right)^{\frac q2}$. However, even if $\psi$ is uniformly bounded, the standard theory of quasilinear equation does not apply since $q>2$. We have
$\psi_t=(v_t-\gl v)e^{-\gl t}$ and $\nabla_g\psi=e^{-\gl t}\nabla_gv$, hence 
$$\left((\psi_t+\gl\psi+2e^{-\gl t})^2+|\nabla_g\psi|^2\right)^{\frac q2}\leq C,
$$
which implies 
\bel{Q6-7}\BA{lll}\dsps
|F^*(t)|\leq Ce^{\left((2-q)+(q-1)\gl\right)t}+C'e^{-\gl t}\leq C''.
\EA\ee
Denoting again by $H$ the subspace of functions in $L^2(S^{N-1})$ with zero average, then 
$$\norm{e^{(\frac{2-N}{2}-\gl)t}S(t)[\phi]}_{L^\infty(S^{N-1})}\leq Ce^{(1-N-\gl)t}\norm{\phi}_{L^\infty(S^{N-1})}\quad\text{for all }\phi\in H\cap L^\infty(S^{N-1}),
$$
and
$$\norm{e^{(\frac{N-2}{2}+\gl)\gt}S(\gt)[\phi]}_{L^\infty(S^{N-1})}\leq Ce^{(\gl-1)\gt}\norm{\phi}_{L^\infty(S^{N-1})}\quad\text{for all }\phi\in H\cap L^\infty(S^{N-1}).
$$
These inequalities imply
$$\BA{lll}\dsps \int_0^te^{(\frac{2-N}{2}-\gl)s}\norm{S(s)\left[\int_0^\infty e^{(\frac{N-2}{2}+\gl)\gt}S(\gt)[F^*(t+\gt-s)]d\gt\right]_{L^\infty(S^{N-1})}ds}\leq C.
\EA$$
Therefore the absolute convergence of the term $\dsps\int_0^\infty e^{(\frac{N-2}{2}+\gl)\gt}S(\gt)[F^*(t+\gt-s)]d\gt$ in $L^\infty(S^{N-1})$ is verified. Finally, we have to check that $\psi^*$ belongs to the uniqueness class of solutions of $(\ref{Q6-5})$ with zero average
$$\xi_{tt}+(N-2+2\gl)\xi_t+\gl(N-2+\gl)\xi+\Gd_g\xi=0
$$
in the space $H$. The associated ODE is 
$$x''+(N-2+2\gl)x'+(1-N+\gl(N-2+\gl))x=0
$$
with characteristic equation 
$$\gr^2+(N-2+2\gl)\gr+1-N+\gl(N-2+\gl)=0
$$
and roots
$$\gr_+=1-\gl\quad\text{and }\gr_-=1-\gl-N.
$$
Hence, if the function $\psi^*$ satisfies $\norm{\psi^*}_{L^2(S^{N-1})}=o(e^{\gr_+t})$, it belongs to the uniqueness class associated to $(\ref{Q6-5})$ and the integral representation $(\ref{Q6-6})$ is valid. The function $\psi$ is bounded and thus the representation holds.
This ends the proof.\qeda

\subsection{Proof of \rth{th7}}

The proof is the consequence of the two next intermediate results

\blemma{th*3} Let $N\geq 3$ and $q>2$. If $u$ is a solution of $(\ref{Na-1})$ in $B_1^c$ such that 
\bel{Q6-8*}
0<m_1\leq |x|^2e^{u(x)}\leq M_1 \quad\text{for all }x\in B_2^c,
\ee
there exists a solution $\gw$ of $(\ref{Na-14})$  such that 
\bel{Q6-9*}
u(r,\gs)+2\ln r\to\gw(\gs)\quad\text{as }|x|\to\infty,
\ee
uniformly on $S^{N-1}$.
\es
\Proof The function $v$ defined in $(\ref{Q6-3*})$ is bounded, but the estimate $(\ref{Q6-4**})$ on $v_t$ and $\nabla'v$ is not enough to use the energy method as in the case $1<q<2$, and the standard regularity theory for quasilinear elliptic equations does not apply directly since $q>2$. \\
\nind{\it Step 1: we claim that $|v_{tt}|+|D^2v|+|v_t| +|\nabla'v|$ is uniformly bounded.}
We write $(\ref{Q6-3})$ under the form 
$$v_{tt}+(N-2)v_t+\Gd_g v-K-D(t,\gs)\left((v_t+2)^2+|\nabla_gv|^2\right)^{\frac 12}=0,
$$
where $K=e^v-2(N-2)$ is uniformly bounded and
$$D(t,\gs)=me^{(2-q)t}\left((v_t+2)^2+|\nabla_gv|^2\right)^{\frac {q-1}2},
$$
verifies $|D(t,\gs)|\leq Me^{(2-q)t}e^{(q-2)t}=M$ independently of $(t,\gs)$ by $(\ref{Q6-4**})$. For this equation the standard regularity theory applies and we get
$$\norm {v(.,.)}_{C^2((T-1,T+1)\ti S^{N-1})}\leq C_1\norm {v(.,.)}_{L^\infty((T-2,T+2)\ti S^{N-1})}+C_2\leq C_3\quad\text{for all }T>2.
$$
\nind{\it Step 2: end of the proof} As in the proof of \rprop{conv} the energy identity is verifies
$$\BA{lll}\dsps\frac{1}{2}\frac{d}{dt}\int_{S^{N-1}}\left(v_t^2-|\nabla_gv|^2-2(N-2)v+e^v\right) dS\\[4mm]
\phantom{----}\dsps=(N-2)\int_{S^{N-1}}v_t^2dS+me^{(q-2)t}\int_{S^{N-1}}\left((v_t+2)^2+|\nabla_gv|^2\right)^{\frac q2}v_tdS.
\EA$$
Hence 
$$\int_2^\infty\int_{S^{N-1}}v_t^2dSdt<\infty.
$$
The trajectory of $v(t,.)$ is relatively compact and we end the proof using \rth{th**} extracted from \cite{HuTak}.\qeda
\blemma{th*4} Let $N\geq 3$ and $q>2$. If $u$ is a solution of $(\ref{Na-1})$ in $B_1^c$ such that $|x|^2e^{u(x)}$ remains bounded in $B_2^c$ then $(\ref{Q6-8*})$ holds for some 
$M_1>m_1>0$.
\es
\Proof The equivalent form of $(\ref{Q6-8*})$ is
\bel{Q6-8}\BA{lll}\dsps
|2\ln |x|+u(x)|\leq C\quad\text{for all }|x|\geq 2.
\EA\ee
{\it Step 1: we claim that}
\bel{Q6-9}\BA{lll}\dsps
 -c^*\ln |x|-C_2\leq u(x)\leq -2\ln |x|+C_1\quad\text{for }|x|\geq 1,
\EA\ee
for some $c^*\geq 2$ and $C_i\in\BBR$.\\
The average $\bar u$ satisfies
$$-\bar u_{rr}-\frac{N-1}{r}\bar u_r\leq \overline {e^u}\leq \frac{M_1}{r^2}.
$$
Therefore 
$$-(r^{N-1}\bar u_r)_r\leq M_1r^{N-3}.
$$
By integration on $(1,r)$
$$-r^{N-1}\bar u_r(r)+\bar u_r(1)\leq \frac{M_1}{N-2}(r^{N-2}-1).
$$
Hence 
$$-\bar u_r(r)\leq \frac{M_1}{(N-2)r}- r^{1-N}\left(\bar u_r(1)+\frac{M_1}{N-2}\right).
$$
Integrating again we obtain
$$-\bar u(r)+\bar u(1)\leq \frac{M_1}{N-2}\ln r+\left(\bar u_r(1)+\frac{M_1}{N-2}\right)\frac{r^{2-N}-1}{N-2}.
$$
Therefore 
\bel{Q6-10}\BA{lll}\dsps
\bar u(r)\geq -\frac{M_1}{N-2}\ln r+C_1,
\EA\ee
which in turn implies by \rprop{th*3}
\bel{Q6-11}\BA{lll}\dsps
u(x)\geq -\frac{M_1}{N-2}\ln |x|-C_2= -c^*\ln |x|-C_2.
\EA\ee
\nind{\it Step 2: we claim that}
\bel{Q6-11*}\BA{lll}\dsps
|v_t|+|\nabla'v|+|v_{tt}|+|D^2v|\leq Ct\quad\text{for }t\geq T,
\EA\ee
where $v(t,\gs)=u(r,\gs)+2t$ is defined in $(\ref{Q6-3*})$. \\
We set $t\phi(t,\gs)=v(t,\gs)$. Since $(2-c^*)t-C_2\leq v(t,\gs)\leq C_1$, the function $\phi$ is bounded in $[T,\infty)\ti S^{N-1}$ where it satisfies
$$\phi_{tt}+\left(N-2+\frac 2t\right)\phi_t+\Gd_g\phi-\frac{2(N-2)}{t}+\frac 1te^{t\phi}+mt^{q-1}e^{(2-q)t}\left(\left(\phi_t+\frac 1t(\phi-2)\right)^2+|\nabla'\phi|^2\right)^{\frac q2}=0.
$$
As a rough upper estimate we have $|\phi_t|+|\nabla'\phi|\leq e^{\frac{q-2}{q-1}t}$. Then we write
$$\left(\left(\phi_t+\frac 1t(\phi-2)\right)^2+|\nabla'\phi|^2\right)^{\frac q2}=\left(\left(\phi_t+\frac 1t(\phi-2)\right)^2+|\nabla'\phi|^2\right)\left(\left(\phi_t+\frac 1t(\phi-2)\right)^2+|\nabla'\phi|^2\right)^{\frac {q-2}2},
$$
and
$$\left(\left(\phi_t+\frac 1t(\phi-2)\right)^2+|\nabla'\phi|^2\right)^{\frac {q-2}2}\leq Ce^{\frac{(q-2)^2}{q-1}t}.
$$
If we set 
$$G(t,.)=mt^{q-1}e^{(2-q)t}\left(\left(\phi_t+\frac 1t(\phi-2)\right)^2+|\nabla'\phi|^2\right)^{\frac {q-2}2},
$$
then $G(t,.)\to 0$ uniformly as $t\to\infty$ and thus
$$0\leq G(t,.)\leq Ct^{q-1}e^{\frac{2-q}{q-1}t}\quad\text{for }t\text{ large enough}.
$$
Now we write the equation satisfied by $\phi$ under the form
\bel{Q6-12}\BA{lll}\dsps
\phi_{tt}+\left(N-2\right)\phi_t+\Gd_g\phi-\frac{2(N-2)}{t}+\frac 1te^{t\phi}+\frac 2t\phi_t+G(t,.)\left(\left(\phi_t+\frac 1t(\phi-2)\right)^2+|\nabla'\phi|^2\right)=0.
\EA\ee
Since $\phi$ is bounded, the classical regularity theory for quasilinear elliptic equations with quadratic growth in the gradient applies and we deduce that 
$|\phi_t|$, $|\nabla'\phi|$ $|\phi_{tt}|$ $|D^2\phi|$ and $|\nabla'\phi_t|$ remain bounded. Since $v_t=t\phi_t+\phi$, $\nabla'v=t\nabla \phi$, etc. we obtain $(\ref{Q6-11})$.\smallskip

\nind{\it Step 3: end of the proof.} Since $v(t,.)-\bar v(t)$ is uniformly bounded, there exists $\gth>0$ such that $\overline {e^v}\leq \gth e^{ \bar v}$, hence 
$$\bar v_{rr}+(N-2)\bar v_{r}-2(N-2)+\gth e^{ \bar v}\geq 0.
$$
Assume first that $\bar v(t)\to-\infty$ when $t\to\infty$. Hence for $t\geq t_0>\ln 2$, we have $2(N-2)-e^{\gth \bar v}>1$ since $N\geq 3$, hence
$$\left(e^{(N-2)t}\bar v_{r}\right)'\geq e^{(N-2)t},
$$
which implies 
$$\bar v_{r}(t)\geq \frac{1}{N-2}+e^{(N-2)(t_0-t)}\left(\bar v_{r}(t_0)-\frac{1}{N-2}\right).
$$
This implies that $\bar v(t)\to\infty$ when $t\to\infty$, contradiction.\\
Assume now that  
$$\liminf_{t\to\infty}\bar v(t)=-\infty\quad\text{and }\limsup_{t\to\infty}\bar v(t)=K>-\infty.
$$
Then there exists a sequence $\{t_n\}$ converging to $\infty$ such that $\bar v(t_{2n})\to -\infty$, $\bar v_{r}(t_{n})=0$, $\bar v(t)$ is nonincreasing on $(t_{2n-1},t_{2n})$
 and $\bar v(t_{2n})<\bar v(t_{2n-1})$. Then, multiplying by $\bar v_t$ and integrating on $(t_{2n-1},t_2n)$ yields
 $$\BA{lll}\dsps\frac{1}{2}\underbrace{\left(\bar v_{r}^2(t_{2n})-\bar v_{r}^2(t_{2n-1})\right)}_{=0}+(N-2)\int_{t_{2n-1}}^{t_{2n}}\bar v_{r}^2(t)dt-2(N-2)\left(\bar v(t_{2n})-\bar v(t_{2n-1})\right)\\[2mm]
 \phantom{--------------------------}+\gth \left(e^{\bar v(t_{2n})}-e^{\bar v(t_{2n-1})}\right)\leq 0.
 \EA$$
 We encounter two possibilities:\\
 (I) $\bar v(t_{2n})-\bar v(t_{2n-1})$ is unbounded and we have a contradiction with the fact that $e^{\bar v(t_{2n})}-e^{\bar v(t_{2n-1})}$ is bounded since $v\leq C_1$.\\
 (II) $\bar v(t_{2n})-\bar v(t_{2n-1})$ is bounded. Then  
 $$e^{\bar v(t_{2n})}-e^{\bar v(t_{2n-1})}=e^{\gt \bar v(t_{2n})+(1-\gt)\bar v(t_{2n-1})}\left(\bar v(t_{2n})-\bar v(t_{2n-1})\right),
$$
and we obtain
$$(N-2)\int_{t_{2n-1}}^{t_{2n}}\bar v_{r}^2(t)dt-\left(2(N-2)-e^{\gt \bar v(t_{2n})+(1-\gt)\bar v(t_{2n-1})}\right)\left(\bar v(t_{2n})-\bar v(t_{2n-1})\right)\leq 0.
$$
Since $2(N-2)-e^{\gt \bar v(t_{2n})+(1-\gt)\bar v(t_{2n-1})}>0$ when $n\to\infty$, we obtain a contradiction. This ends the proof.\qeda
  
\begin{center}

{\bf STATEMENTS AND DECLARATION}
\end{center}

\nind Conflict of interest: The authors declare that they have no conflict of interest and that this
article follows all the ethical rules.\\

\nind Data availability: Data sharing is not applicable to this article as no datasets were generated
or analysed during the study.

\nind{\bf Aknowledgements}
The authors are grateful to the referee for careful checking of the manuscript and
relevant observations and suggestions.\medskip

\nind Marie-Fran\c{c}oise Bidaut-V\'eron \\
Institut Denis Poisson CNRS UMR 7013. \\
Universit\'e de Tours\\
Tours, France\\
{\it veronmf@univ-tours.fr}\\

\nind Laurent V\'eron \\
Institut Denis Poisson CNRS UMR 7013\\
Universit\'e de Tours\\
Tours, France\\
{\it veronl@univ-tours.fr}\\
{\it laurent.veron49@gmail.com}

 \end{document}